\documentclass[12pt,a4paper]{article}

\usepackage{a4,enumerate}
\usepackage{amsthm}
\usepackage{amssymb}
\usepackage{amsmath}
\usepackage{latexsym}
\usepackage{pstricks,pst-coil,pst-plot,pst-node}
\usepackage{graphicx}

\long\def\bild#1{#1}
\def\url#1{{\tt #1}}

\newdimen\templaenge

\DeclareMathAlphabet{\doba}{U}{msb}{m}{n}
\gdef\mC{\doba{C}}

\gdef\mN{\doba{N}}

\gdef\mR{\doba{R}}
\gdef\mS{\doba{S}}

\gdef\mZ{\doba{Z}}

\def\qed{{\leavevmode\unskip\nobreak\hfil\penalty 50\hskip 1em%
  \hbox{}\nobreak\hfil\lower 1pt\hbox{$\Box$\kern-.5pt}\parfillskip 0pt
  \finalhyphendemerits 0\par\bigbreak}}
\def\qedmath#1{\setbox0\hbox{$\displaystyle #1$}\templaenge=\textwidth\advance\templaenge by -\wd0%
\setbox1\hbox{$\Box$}\advance\templaenge by -2\wd1%
$$#1\hbox to0pt{\kern.5\templaenge$\Box$\kern-.5pt\hss}$$\par\bigbreak}

\def\al{{\alpha}}
\def\be{{\beta}}
\def\de{{\delta}}

\def\om{{\omega}}

\def\la{{\lambda}}

\def\si{{\sigma}}
\def\Si{{\Sigma}}
\def\ga{{\gamma}}
\def\ep{{\varepsilon}}
\def\Ga{{\Gamma}}

\def\ph{{\varphi}}
\def\phi{{\varphi}}

\def\na{{\nabla}}
\def\pa{{\partial}}

\def\cS{\mathcal{S}}

\def\ohne{\smallsetminus}
\def\ti{\tilde}
\def\witi{\widetilde}

\def\ol{\overline}

\def\embed{\hookrightarrow}

\def\Dflat{D^{\rm eucl}}

\def\Ree{{\mathop{\rm Re\;}}}
\def\Imm{{\mathop{\rm Im\;}}}

\def\ie{i.\thinspace e.\ \ignorespaces}

\def\nummerarray#1#2{\par\noindent\setbox0\hbox{\rm (#1)}\setbox1\hbox{$#2$}\unhcopy0%
\dimen0=.5\textwidth \advance\dimen0 by -\wd0 \advance\dimen0 by -.5\wd1 \kern\dimen0 \unhcopy1}

\def\ker{\mathop{{\rm ker}}}

\def\dim{\mathop{{\rm dim}}}
\def\dist{\mathop{{\rm dist}}}
\def\vol{{\mathop{{\rm vol}}}}
\def\dvol{{\mathop{{\rm dvol}}}}
\def\area{{\mathop{{\rm area}}}}
\def\diam{{\mathop{{\rm diam}}}}

\def\min{\mathop{{\rm min}}}
\def\max{\mathop{{\rm max}}}

\def\mod{\mathop{{\rm mod}}}
\def\can{\mathop{{\rm can}}}

\def\Id{\mathop{{\rm Id}}}
\def\End{{\mathop{{\rm End}}}}

\def\SU{{\rm SU}}

\def\Hom{{\mathop{{\rm Hom}}}}
\def\scal{{\mathop{{\rm scal}}}}

\def\geucl{{g_{\mathop{\rm eucl}}}}

\def\res#1#2{{#1}\lower .11ex\hbox{$|$}\lower .644ex\hbox{$\scriptstyle #2$}}

\def\lammin{{\lambda_{\rm min}^+}}

\long\def\komment#1{}
\long\def\longver#1{}

\def\proof#1{{\par\medbreak\noindent {\bf Proof\setbox0\hbox{#1}%
\ifdim\wd0=0pt .\else\ \ignorespaces #1.\fi}\enspace}}
\def\iop#1{{\par\medbreak\noindent {\bf Idea of proof\setbox0\hbox{#1}%
\ifdim\wd0=0pt .\else\ \ignorespaces #1.\fi}\enspace}}
\def\examples{{\noindent {\bf Examples. }\par\kern-\baselineskip}}

\newtheoremstyle{remarks}{3pt}{3pt}{}{}{\bfseries}{}{ }{}
\newtheoremstyle{appendix}{3pt}{3pt}{}{}{\bfseries}{}{ }{}

\makeatletter
\@addtoreset{equation}{section}
\renewcommand{\theequation}{\the\c@section.\the\c@equation}
\makeatother
\newtheorem{theorem}[equation]{\bf T{\footnotesize HEOREM}}
\newtheorem{proposition}[equation]{\bf P{\footnotesize ROPOSITION}}
\newtheorem{lemma}[equation]{\bf L{\footnotesize EMMA}}

\newtheorem{conjecture}[equation]{\bf C{\footnotesize ONJECTURE}}
\newtheorem*{cmcprinciple}{\bf Principle for construction of cmc-surfaces}

\newtheorem{theorema}{\bf T{\footnotesize HEOREM}}[section]
\newtheorem{lemmaa}[theorema]{\bf L{\footnotesize EMMA}}
\newtheorem{propositiona}[theorema]{\bf P{\footnotesize ROPOSITION}}
\theoremstyle{definition}
\newtheorem*{remark}{Remark}

\newtheorem*{definition}{Definition}

\theoremstyle{remarks}

\def\eref#1{{\rm (\ref{#1})}}

\def\pD{{p_D}}
\def\pY{{p_Y}}
\def\qD{{q_D}}
\def\Dir{D}

\def\cF{{\cal F}}
\let\<\langle
\let\>\rangle
\let\La\Lambda
\def\inj{\mathop{\rm inj}}

\def\extconf{\ol{[g_0]}}

\def\RHS{{\rm RHS}}

\begin{document}

\title{The smallest Dirac eigenvalue in a spin-conformal class
and cmc-immersions}
\author{Bernd Ammann%
}
\date{September 2005}
\maketitle

\begin{abstract}
Let us fix a conformal class $[g_0]$ and a spin structure $\si$ on
a compact manifold $M$. For any $g\in [g_0]$, let $\la^+_1(g)$
be the smallest positive eigenvalue of the Dirac operator $D$ on $(M,g,\si)$.
In a previous article we have shown that
  $$\lammin(M,g_0,\si):=\inf_{g\in [g_0]} \la_1^+(g)\vol(M,g)^{1/n}>0.$$
In the present article, we enlarge the conformal class by adding 
certain singular metrics. We will show that
if $\lammin(M,g_0,\si)<\lammin(S^n)$,
then the infimum is attained on the enlarged conformal class.
For proving this, we solve a system of
semi-linear partial differential equations involving a nonlinearity
with critical exponent:
  $$D\phi= \la |\phi|^{2/(n-1)}\phi.$$
The solution of this problem has many analogies
to the solution of the Yamabe problem.
However, our reasoning is more involved than in the Yamabe problem
as the eigenvalues of the Dirac operator tend to $+\infty$ and $-\infty$.

Using the spinorial Weierstra\ss{} representation, 
the solution of this equation in dimension 2 shows the existence 
of many periodic constant mean curvature surfaces.

\end{abstract}

{\bf Mathematics Classification:} 58J50, 53C27 (Primary),
58C40, 35P15, 35P30, 35B33 (Secondary)

\section{Introduction}



Let $M$ be a compact $n$-dimensional 
manifold, $n\geq 2$, with a fixed conformal class $[g_0]$ and a fixed 
spin structure $\si$.
Let $g$ be conformal to $g_0$, i.e.\ $ g\in [g_0]$. 
The classical Dirac operator $D_g$
on $(M,g,\si)$ has discrete real spectrum
with finite multiplicities. The eigenvalues tend to $+\infty$ and $-\infty$.
The dimension of the kernel of $D_g$ is a spin-conformal invariant, i.e.
$\dim \ker D_g=\dim \ker D_{g_0}$ for all $g\in [g_0]$.
We denote the first (=smallest) positive eigenvalue of $D_g$ 
by $\la_1^+(g)$,
and the first (=largest) negative eigenvalue by $\la_1^-(g)$.

Finding bounds for this eigenvalue has attracted much interest
during the last decades. Among several estimates 
\cite{lichnerowicz:63,friedrich:80,kirchberg:86,kirchberg:88,
kramer.semmelmann.weingart:98,kramer.semmelmann.weingart:99} 
in terms of a positive scalar curvature bound, 
let us mention the following estimate due to Friedrich~\cite{friedrich:80}.
If the minimum of the scalar curvature
of $(M,g)$ is at least $s>0$, then 
$\left(\la_1^\pm(g)\right)^2\geq {n\over 4(n-1)}s$.
  
An improvement of this inequality that is important in conformal 
geometry was derived by Hijazi \cite{hijazi:86} 
in terms of the first eigenvalue
$\la_1(L_g)$ 
of the conformal Laplacian $L_g=4\,{n-1\over n-2}+\scal_g$, 
namely
\begin{equation}\label{ineq.hij}
  \left(\la_1^\pm(g)\right)\geq {n\over 4(n-1)} \la_1(L_g)
\end{equation}
if $n\geq 3$. On the other hand, if the 
\emph{Yamabe invariant}
 $$Y(M,[g_0])=\inf_{g\in [g_0]}{\int_M \scal_g\dvol_g\over
\vol(M,g)^{(n-2)/n}}\in (-\infty, n(n-1)]$$ 
is non-negative, one can use the  conformal
transformation formula
for the scalar curvature \cite[Theorem~1.159]{besse:87}, and one obtains
that 
  $$ \la_1(L_g)\vol(M,g)^{2/n}\geq Y(M,[g_0]).$$
Together with \eref{ineq.hij}, we get
\begin{equation}\label{ineq.hij.der}
\left|\la_1^\pm(g)\right|\vol(M,g)^{1/n} \geq 
    \sqrt{{n\over 4(n-1)}\, Y(M,[g_0])}.
\end{equation}
That this inequality even holds in 
dimension $n=2$, was proved by C.~B\"ar \cite{baer:92b}. In this special case,
the Gauss-Bonnet theorem tells us that $Y(M)= 4\pi \chi(M)$, hence we obtain
\begin{equation}\label{ineq.baer}
 \left|\la_1^\pm(g)\right|\area(M,g)^{1/2}\geq 2\sqrt{\pi}
\end{equation}
if $M$ is diffeomorphic to $S^2$, but we do not get a bound for other surfaces.
Equality in \eref{ineq.hij.der} and \eref{ineq.baer} hold for the round 
spheres.

Hence, we have obtained an explicit positive 
lower bound $\left|\la_1^\pm(g)\right|\vol(M,g)^{1/n}$ that is uniform
on the conformal class $[g_0]$.
From now on, we will restrict to the first \emph{positive} eigenvalue, 
in order to simplify the presentation (see the remark below).

The \emph{existence} of a positive lower 
bound for $\la_1^+(g)\vol(M,g)^{1/n}$
had already been derived before by J. Lott. His estimate
does not require that $Y(M,[g_0])$ is positive, but needs the weaker assumption
that the Dirac operator is invertible, i.e.\ has $0$-dimensional kernel.
Unfortunately, his bound is not explicit, and for most spin-conformal
manifold the determination of the 
value of the associated infimum 
\begin{equation}\label{def.lammin}
\lammin(M,[g_0],\si)
    :=\inf_{g\in [g_0]} \la_1^+(g)\vol(M,g)^{1/n}
\end{equation}
is still a challenging open problem. 

A key idea of J. Lott's article is to derive a lower 
bound for the first Dirac eigenvalue in terms of the  
supremum of a conformally invariant functional, a version of this 
functional will be explained in Section~\ref{sec.varprin}.  

In our article \cite{ammann:03}, we started to study this
functional in more detail. In particular, we showed that Lott's result
extends to the case that the Dirac operator has non-trivial kernel,
namely 
  $$\lammin(M,[g_0],\si)>0$$
for any compact Riemannian spin manifold $(M,[g_0],\si)$.

Furthermore, it was  shown that
  $$\lammin(M,[g_0],\si) \leq \lammin(\mS^n),$$
where $\mS^n=(S^n,g_{\rm can})$ denotes the sphere with its standard 
Riemannian metric of constant sectional curvature $1$.
This bound has been proven in \cite{ammann:03} unless 
$\ker D\not=\{0\}$ and $n=2$. The remaining case $\ker D\not=\{0\}$ and $n=2$ 
was given in \cite{grosjean.humbert:p05}.

In the present article, we discuss whether the associated infimum in 
\eref{def.lammin} is attained. For having a well-behaved
minimization problem, it is reasonable to replace the conformal class
$[g_0]$ in 
\eref{def.lammin} by its $L^\infty$-completion $\overline{[g_0]}$.
Here, by definition a metric $f^2g_0$ is in $\overline{[g_0]}$ if 
$f$ is a real-valued $L^\infty$ function. The first positive 
eigenvalue of the Dirac operator extends naturally to this completion
(see Section~\ref{sec.spd}).
 We show that if we have the strict inequality
 \begin{equation}\label{eq.strict}
  \lammin(M,g_0,\si)< \lammin(\mS^n)= {n\over 2}\,\om_n^{1/n},
\end{equation}
then the infimum is obtained by a generalized metric $g\in \overline{[g_0]}$.
This minimizer has the form $g:=|\phi|^{4/(n-1)}g_0$ where 
$\phi$ is a spinor of regularity $C^2$. The set 
$\phi^{-1}(0):=\{x\in M\,|\,\phi(x)=0\}$ is  
called the \emph{nodal set of $\phi$} or the 
\emph{set of degeneration of $g$}.

\begin{theorem}\label{theo.main}
Let $M$ be a compact manifold of dimension $n\geq 2$
with a fixed conformal
class $[g_0]$ and a spin structure $\si$.
Assume that \eref{eq.strict}
holds. Let $\al:= 2/(n-1)$ if $n\geq 4$, and let $\al\in (0,1)$ if $n\in\{2,3\}$.  
\begin{enumerate}[\rm (A)]
\item Then there is a spinor field
   $\phi\in C^{2,\al}(\Si M)\cap C^\infty(\Si(M\ohne \phi^{-1}(0)))$ 
   on $(M,g_0)$ such that
\begin{equation}\label{eq.dirac.nonlin}
  D_{g_0}\phi = \lammin\, |\phi|^{2/(n-1)}\phi,\qquad \|\phi\|_{2n/(n-1)}=1.
\end{equation}
\item There is a $g\in \extconf$ such that
  $$\la^+_1(g)\vol(M,g)^{1/n}=\lammin.$$
The metric has the form $g=|\phi|^{4/(n-1)}g_0$ where 
$\phi$ is a spinor as in (A).
\item If $\dim M=2$, then the metric $g$ is smooth 
and the set of degeneration of $g$, denoted $\cS_g$,
is finite.
Furthermore
  $$\#\cS_g< {\rm genus}(M).$$
In particular,  if $M$ is diffeomorphic to a
$2$-torus, then the set of degeneration $\cS_g$ is empty.
\end{enumerate}
\end{theorem}

Inequality \eref{eq.strict} has been proven for several 
classes of manifolds.
It is known that it holds for non-conformally-flat manifolds of dimension
$\geq 7$ \cite{ammann.humbert.morel:p03av2}. 

Inequality~\eref{eq.strict} has also been shown if 
$M$ is conformally flat, if $D$ is invertible, 
and if the mass endomorphism is  not identically zero 
(after a possible change of orientation if $\dim M \equiv 3 \mod 4$) 
\cite{ammann.humbert.morel:p03b}. 
The mass endomorphism is a section of $\End(\Si M)$
defined as the zero order term of the development of the Green function
for the Dirac operator at the diagonal with respect to a 
conformal coordinate map.

Furthermore, \eref{eq.strict} is known for many Riemann surfaces (i.e. $n=2$), 
e.g.\ all rectangular tori have a spin structure such that 
\eref{eq.strict} holds.

The Euler-Lagrange equation~\eref{eq.dirac.nonlin} of the above minimization
problem has a particularly nice  interpretation in dimension~$n=2$. 
Locally the equation~\eref{eq.dirac.nonlin} 
can be translated into a conformal constant mean curvature
immersion into $\mR^3$. 
This translation
is a spinorial extension of the Weierstrass representation (see
Section~\ref{sec.weierrep}). By pasting together these local surfaces
one obtains a ``periodic branched conformal cmc immersion''. 

More exactly, let $(M,g)$ be a compact Riemann surface together 
with its universal covering $\pi:\witi M\to M$. A 
\emph{periodic branched conformal cmc immersion based on $(M,g)$}
is by definition a smooth map $F:\witi M\to \mR^3$ together with 
finitely many points $p_1,\ldots,p_k\in M$, the so-called 
\emph{branching points}, such 
that the following properties hold:
\begin{enumerate}[{\rm (1)}]
\item {\it Periodicity:}  There is a homomorphism 
$h:\pi(M)\to \mR^3$, the \emph{periodicity map}, such that 
for any $\gamma\in \pi_1(M)$, and $x\in\witi M$ one has
  $$F(x\cdot \gamma)=F(x)+h(\gamma).$$
Here $\cdot$ denotes the action of $\pi_1$ on $\witi M$ via Deck transformation.
\item {\it Conformality:} The restriction of $F$ to $\witi M\setminus \pi^{-1}(\{p_1,\ldots,p_k\})$
      is a conformal immersion.
\item {\it Branching points:} We have 
$dF_q=0$ for any $q\in \pi^{-1}(\{p_1,\ldots,p_k\})$. 
The order of the first non-vanishing term in the 
Taylor development of $dF$ in $q$ is called the 
\emph{branching index} of $F$ at $q$. 
\item {\it CMC:} The image $F\left(\witi M\setminus \pi^{-1}(\{p_1,\ldots,p_k\}\right)$ 
is an immersed surface 
with constant mean curvature.
\end{enumerate}

The principle yields the existence of many periodic branched conformal
cmc immersions. The set of periodic branched conformal
cmc immersions, $H\neq 0$, is 
essentially (see Section~\ref{sec.app.cmc} for a precise 
statement) 
in bijection with the stationary points of the variational problem
associated to our minimization problem. 
In particular, all minimizers of \eref{def.lammin} give rise
to periodic branched conformal
cmc immersions. We obtain

\begin{cmcprinciple}
Assume that the Riemann spin surface $(M,g,\si)$ carries a metric
$g$ such that the first positive eigenvalue of the Dirac operator is smaller
than  $2\sqrt{\pi/\area(M,g)}$. Then there is a periodic branched
conformal cmc immersion $F$ based on $(M,g)$. 
The regular homotopy class of $F$ is determined by the spin-structure $\si$.
The indices of all branching points are even, 
and the sum of these indices
is smaller than $2 {\rm genus}(M)$. In particular, 
if $M$ is a torus, there are no branching points.  
\end{cmcprinciple}

Some examples of branched
conformal cmc immersions that may arise by this principle
are given in Section~\ref{sec.app.cmc}. 

The problem that we are discussing in this article has many relations to the 
Yamabe problem. For a given compact conformal manifold $(M,[g_0])$ 
the Yamabe problem is the problem to find a metric of constant scalar curvature
$g$ in $[g_0]$. The problem has been affirmatively solved by Trudinger, 
Aubin, Schoen and Yau, see \cite{lee.parker:87} for a good overview. 
At first, due to the 
conformally invariant character of our problem and the Yamabe problem,
certain bounded, but non-compact Sobolev embeddings play an important role.
In the solution of both problems it is useful to break the conformal 
invariance by perturbing a parameter. For this perturbed parameter
the embeddings are compact, and standard methods yield the existence 
of a minimizer. It then has to be checked whether the perturbed minimizers
converge to a minimizer of the unperturbed problem, or whether they concentrate
in some points. In the Yamabe problem
the minimizers converge if 
\begin{equation}\label{eq.yamab}
Y(M,[g_0])<Y(\mS^n),
\end{equation}
holds whereas in our problem, the minimizers converge if 
\eref{eq.strict} holds.
Secondly, \eref{ineq.hij.der} implies that any manifold satisfying
\eref{eq.strict} satisfies \eref{eq.yamab} as well. Hence, proving 
\eref{eq.strict} for a given spin-conformal manifold $(M,[g_0],\si)$ 
solves the Yamabe problem on this manifold.
The third relationship is that some proofs of inequality \eref{eq.strict}
in special cases resemble to proofs of inequality, 
e.g.\ \cite{ammann.humbert.morel:p03av2} resembles to \cite{aubin:76}.

The structure of the article is as follows.
In Section~\ref{sec.varprin} we reformulate our problem 
as a variational problem. We will see that it is natural to admit 
in the infimum \eref{def.lammin} certain singular metrics namely ``metrics''
conformal to $g_0$ whose conformal factor might vanish 
somewhere. Such metrics --- called \emph{generalized metrics} ---
are the subject of Section~\ref{sec.spd}.  In the following section, 
we discuss the round sphere. This example is helpful to obtain a deeper 
understanding for the analytical difficulties. However, it can be skipped
if the reader is only interested in the main results of the article.
Section~\ref{sec.reg.theo} is devoted to certain regularity issues that will 
become important in Section~\ref{sec.solution}, where 
the variational principle is finally solved. We then show in 
Section~\ref{sec.proofmain} how this implies the main theorem.
The singularities of the minimizers are discussed in Section~\ref{sec.degen}.
In Section~\ref{sec.weierrep} we recall the spinorial 
Weierstrass representation. Section~\ref{sec.app.cmc} uses the spinorial 
Weierstrass representation to derive the application to constant mean 
curvature surfaces. Several examples are included.
In Appendix~\ref{sec.regularity} we summarize (without proofs)
some analytical tools. 
This appendix shall also serve as a reference for fixing the notations for 
Sobolev spaces and H\"older spaces. In Appendix~\ref{sec.scha.power} 
we proof a proposition about H\"older spaces, as we could not find a 
proof in the literature.

\begin{remark}[The first negative eigenvalue]
The infimum is also attained in \eref{def.lammin}
if we replace $\la_1^+$ by $|\la_1^-|$.  
In the case $n\not\equiv 3\mod 4$ this is obvious: 
there exists an automorphism
of the spinor bundle anticommuting with $D$, 
hence $\la_1^-(g)=-\la_1^+(g)$. 
In the case $n\equiv 3\mod 4$, the proof for $\la_1^-$
is up to the obvious sign changes completely identical.
In almost all references cited in the introduction, 
all statements for   $\la_1^+$ also hold for $|\la_1^-|$ and vice versa.
The only exception is \cite{ammann.humbert.morel:p03b}.
\end{remark}

{\bf Acknowledgements.}
I am very much indebted to Christian B\"ar for various support.
Thank you also to Emmanuel Humbert for his continuing interest in conformal 
spin geometry and his deep analytical ideas. Many electronic and personal 
discussions with Robert Kusner and Karsten Grosse-Brauckmann were also very helpful. I also want to thank B.~Booss-Bavnbek, Oussama Hijazi,
Sergiu Moroianu,Victor Nistor, Reiner Sch\"atzle and Guofang Wang
for several stimulating discussions related to this article.

\section{The associated variational principle}\label{sec.varprin}

The goal of this section is to reformulate the problem of minimizing
the first Dirac eigenvalue in a conformal class as a variational problem.
The choice of a good functional is not very easy, as 
we would like to find a functional which is both bounded and 
conformally invariant. 
The Dirac operator $D$ has a simple behavior under
conformal change of the metric, its square $D^2$ transforms in a more 
complicated way. Hence, it is desirable to use a functional that 
contains only terms in $D\phi$ and $\phi$ and no term depending
on $D^2\phi$. The $\cF_q$-functional defined in \eref{def.F}, $q=2n/(n+1)$,
is such a conformally invariant, bounded functional, and it will 
be turn out that by working with this functional we obtain
the desired results.

\begin{remark} In case, that the reader of our article
is already familiar with the analytical problems of the Yamabe problem,
he might find it enlightening to 
compare this problem to our problem.
Many techniques from the resolution of the Yamabe problem
can be carried over to our setting. However, several arguments from
the resolution of the Yamabe problem fail. Two main problems arise:
on the one hand the spectrum of $D$ is neither bounded from 
below nor from above, on the other hand, we are working with sections
of a vector bundle instead of functions, hence the standard maximum principle
is not available.
Such arguments will have to be replaced by other approaches.
\end{remark}



At first, we recall some basic definitions and facts from spin geometry.
For details we refer to textbooks as for example 
\cite{lawson.michelsohn:89, roe:88,friedrich:00} 
or to the beautifully written 
self-contained introduction \cite{hijazi:01}.

Let $M$ be a compact manifold equipped with a Riemannian metric
$g_0$ and a spin structure $\si$. Assume that $g$
is a metric conformal to $g_0$. One associates to $(M,g,\si)$ 
a natural complex vector bundle over $M$ called the spinor bundle
$\Si(M,g,\si)\to M$. 
Sections of this bundle are called \emph{spinor fields}
or simply \emph{spinors}. The bundle carries a hermitian metric, a metric connection,
and a Clifford multiplication. These additional structures are used to define
the Dirac operator $D_g:\Gamma(\Si(M,g,\si)\to\Gamma(\Si(M,g,\si)$, which
is a first order elliptic differential operator. 
The Dirac operator $D_g$ 
is essentially self-adjoint, and hence it has a self-adjoint
extension. As a consequence of standard elliptic theory,
the spectrum is real and discrete, and all multiplicities
are finite.

In the following, the spin structure $\si$ will be fixed (and often suppressed
in the notation), whereas the metric $g$ varies inside the conformal class
$[g_0]$. Some objects 
will be noted with an index $g$, which means that they are
defined with respect to~$g$ whereas the same object without the index
$g$ indicates that it is defined with respect to the fixed background metric 
$g_0$ (with some exceptions that are apparent from the context). 
For example, $\dvol_g$ is the volume element associated to $g$ 
and $\dvol=\dvol_{g_0}$ is the one associated to $g_0$.
We will frequently use several norms, as for example
$\|\phi\|_{L^p}$ and $\|\phi\|_{H_1^q}$ 
that are defined with respect to $g_0$ 
unless otherwise stated. We summarize their definition and some important 
analytical tools in Appendix~\ref{sec.regularity}.

The spectrum of the Dirac operator $D_g$ will be denoted as 
    $$\ldots \leq \la_2^-(g)\leq \la_1^-(g)<0=\ldots=0
      <\la_1^+(g)\leq\la_2^+(g)\leq \ldots,$$ 
where each eigenvalue appears with its multiplicity. Note that $0$ 
may be an eigenvalue or not. In the former case $D_g$ is not invertible,
otherwise it is. Elliptic theory shows
$\lim_{k\to\infty}\la_k^+(g) =-\lim_{k\to\infty}\la_{-k}^+(g)=\infty$. 
 
The following transformation formula  will be of central importance.
To our knowledge, the formula was found by Hitchin \cite{hitchin:74}. 
Another reference written up in a more self-contained manner
is \cite{hijazi:01}.

\begin{proposition}[Conformal transformation formula for $D$]\label{prop.conf.change}
Let $g=f^2g_0$, $f:M\to \mR$ smooth and positive.
There is an isomorphism of vector bundles
${F:\Si(M,g_0,\si)\to \Si(M,g,\si)}$
which is a fiberwise isometry such that
  $$\Dir_g(F(\phi))=F\left(f^{-{n+1\over 2}}\Dir_{g_0} f^{{n-1\over 2}}\phi\right).$$
\end{proposition}

It is convenient to define 
  $$\witi{F}(\phi)=F(f^{-{n-1\over 2}}\phi).$$
We will use this isomorphism $\witi{F}$
to identify spinors associated 
to conformal metrics. With this identification the conformal transformation 
formula reads as
  $$\Dir_g(\phi)= f^{-1}\Dir_{g_0}(\phi),$$
and $|\phi|_g=f^{-{n-1\over 2}}|\phi|_{g_0}$.
In particular, with this identification the kernels of $D_g$ and 
$\Dir_{g_0}$ coincide.
It is easy to verify that with this identification,
the expression $\int \<\Dir_g \phi,\phi\>_g\dvol_g$ is conformally invariant. 
The $L^p$-norm $\|\phi\|_{L^p(g)}:=\left(\int |\phi|_g^p\dvol_g\right)^{1/p}$ 
is conformally invariant 
if and only if $p$ has the value $\pD:=2n/(n-1)$. 
Similarly, $\|\Dir_g\phi\|_{L^q(g)}$ is conformally invariant if and only if
$q$ has the value $\qD:=2n/(n+1)$. Note that $\qD^{-1}+\pD^{-1}=1$. 

For any $q\in [\qD,2]$ and for any $H_1^q$-spinor $\phi$ 
that is not in the kernel of the 
Dirac operator we define
\begin{equation}\label{def.F}
  \cF_q^g(\phi)=  {\int\<\Dir_g\phi ,\phi\>_g\dvol_g\over \|\Dir_g\phi\|_{L^q(g)}^2}, \qquad
     \mu_q^g :=\mu_q(M,g,\si):= \sup \cF_q^g(\psi),
\end{equation}
and $\cF_q:= \cF_q^{g_0}$, $\mu_q:=\mu_q^{g_0}$.
The well-definedness of $\cF_q^g$ and some basic properties are given by the 
following lemma.

\begin{lemma}\label{lem.f.elem}
Let $q \in [\qD,\infty)$. Let $\phi$ be a spinor field
of regularity $H_1^q$, $D\phi\neq 0$.
Then,
$\cF_q(\phi)$ is well-defined and real. Furthermore
$\cF_q:H_1^q\setminus \ker D \to\mR$ 
is Fr\'echet differentiable with derivation given by 
\begin{equation}\label{eq.abl}
d\cF_q(\phi)(\psi) = {2\over \|\Dir\phi\|_{L^q}^2}
       \,\int\< \phi  -\rho_{q,\phi} |\Dir\phi|^{q-2}\Dir\phi,\Dir\psi\>,
\end{equation}
where $\rho_{q,\phi}=\cF_q(\phi) \|\Dir\phi\|_{L^q}^{2-q}$.
The supremum $\mu_q$ is positive and finite.
\end{lemma}

From the above considerations it is evident that the 
functional $\cF_q$ is conformally invariant if and only if $q=\qD$.

\proof{}
Let $\phi$ be a spinor field of regularity $H_1^q$, $D\phi\neq 0$. 
Take $p$ with $p^{-1}+q^{-1}=1$, $q\ge {2n\over n+1}$.
Because of $\phi\in H_1^q\embed L^p$, we see with H\"older's inequality
that $\<D\phi,\phi\>$ is integrable. Thus, the numerator of $\cF_q$ is 
well-defined, and hence $\cF_q$ is well-defined. The self-adjointness
of $\Dir$ implies that $\cF_q(\phi)$ is real.
Moreover, because $H_1^q\embed L^p$ is bounded we see that 
$|\int_M \<D\phi,\phi\>|\leq \|D\phi\|_{L^q}\,\|\phi\|_{L^p}$
is bounded from above by a multiple of $\|\phi\|_{H_1^q}^2$. 
Using Theorem~\ref{theo.global.lp} we obtain
$\|\phi\|_{H_1^q}\leq \|D\phi\|_{L^q}+\|\phi\|_{L^q}$, 
and we see that $\cF_q$ is bounded on 
$H_1^q\cap (\ker D)^\perp\setminus \{0\}$. An arbitrary spinor field 
$\phi\in H_1^q$ is written as the sum of $\phi_1\in \ker D$ and a non-zero
$\phi_2\perp\ker D$. As $\cF_q(\phi)=\cF_q(\phi_2)$, we see that
$\mu_q$ is finite. 

Because of $\phi\in H_1^q$, we have 
$|\Dir \phi|^{q-2}\Dir\phi \in L^p$.
Hence the right hand side of \eref{eq.abl} defines a continuous
functional on $L^p$ and thus on $H_1^q$. We denote the functional
by $\psi\mapsto \RHS_\phi(\psi)$.
Similarly, one sees that
  $$\cF_q(\phi+ \psi)-\cF_q(\phi)- \RHS_\phi(\psi)\leq o(\|\psi\|_{H_1^q}),$$
hence $\cF_q$ is Fr\'echet differential with derivative 
$\psi\mapsto \RHS_\phi(\psi)$.

If $\phi$ is an eigenspinor to a positive eigenvalue, then 
$\mu_q\geq \cF_q(\phi)>0$.
\qed

\begin{proposition}[Properties of $\mu_q^g$]\label{prop.func}
The function $\mu_q^g:[\qD,\infty)\to (0,\infty)$ is continuous from the right,
and 
  $$\mu_2(M,g,\si)=(\la^+_1(g))^{-1}$$
Furthermore, if $\vol(M,g)=1$, then $\mu_q^g$ is non-increasing in $q$.
\end{proposition}

\proof{}
We assume $\vol(M,g)=1$, the statements for  $\vol(M,g)\not=1$ 
then follow by rescaling.

That $\mu_q^g$ is non-increasing follows easily from the H\"older inequality. 

In order to show
the continuity from the right, let $q\geq \qD$ be given.
We take a smooth spinor field $\phi$ such that
$\cF_q^g(\phi)\geq \mu_q^g-\ep$.
Observe that  $$\cF_{q'}^g(\phi)= { \|\Dir_g\phi\|_{L^q(g)}^2 \over \|D_g\phi\|_{L^{q'}(g)}^2}\,\cF_{q}^g(\phi).$$
The function $q'\mapsto \|\Dir_g\phi\|_{L^{q'}(g)}$ is continuous, hence
if $q'\geq q$ is sufficiently close to~$q$, then
  $$\mu_{q'}^g\geq \cF_{q'}^g(\psi)\geq \cF_q^g(\psi)-\ep \geq \mu_q^g -2 \ep.$$
Because $q\mapsto \mu_q^g$ is non-increasing, the continuity from the right follows.

The formula 
$\mu_2^g=(\la^+_1(g))^{-1}$ follows directly if one writes~$\phi$ as a sum of eigenspinors and evaluates
$\cF_2^g(\phi)$.
\qed

\begin{figure}
%
%
%
%
\newdimen\axdim
\axdim=1pt
\newdimen\cudim
\cudim=2pt
\def\ticklen{.2}
\def\abst{.5}

\begin{center}
\psset{unit=1cm}

\begin{pspicture}(-1,-\ticklen)(10.5,6)
\psset{linewidth=\axdim}
\psaxes[linewidth=\axdim,labels=none,ticks=none]{->}(0,0)(9.5,5)
\rput(9.9,0){$q$}
\rput(0,5.4){$\mu_q$}
\psline(2.4,-\ticklen)(2.4,\ticklen)
\psline(3.5,-\ticklen)(3.5,\ticklen)
\rput[t](2.4,-\abst){$\qD$}
\rput[t](3.5,-\abst){$2$}
\psline(-\ticklen,2.8)(\ticklen,2.8)
\psline(-\ticklen,3.7)(\ticklen,3.7)
\rput[r](-\abst,2.8){$(\la^+_1)^{-1}$}
\rput[r](-\abst,3.7){$\mu_\qD$}
\psline[linestyle=dotted]{-}(2.4,0)(2.4,3.7)
\psline[linestyle=dotted]{-}(3.5,0)(3.5,2.8)
\psline[linestyle=dotted]{-}(0,3.7)(2.4,3.7)
\psline[linestyle=dotted]{-}(0,2.8)(3.5,2.8)
\psset{linewidth=\cudim}
\psecurve[showpoints=false](2.2,3.74)(2.4,3.7)(3.5,2.8)(4,2.0)(5.1,1.6)(7.2,1.4)(7.5,1.2)
\end{pspicture}
\end{center}
\caption{$\mu_q$ as a function of $q$}
\end{figure}
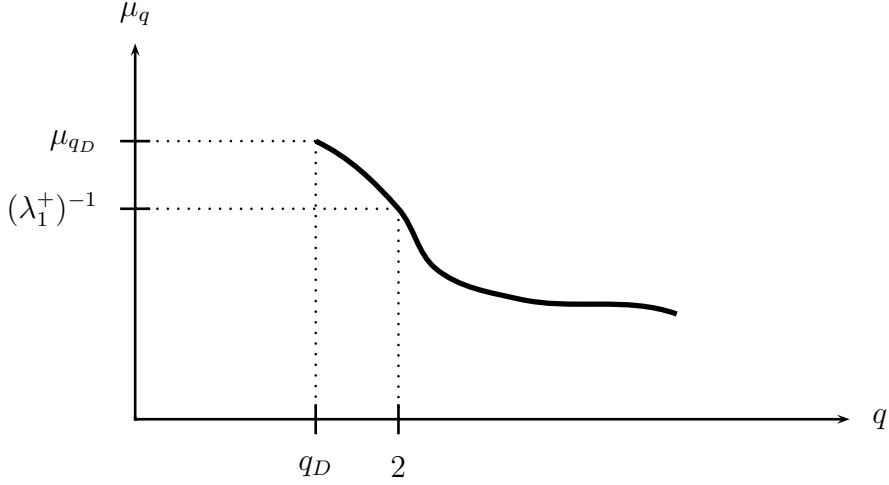

\begin{remark}
In Proposition~\ref{prop.sol.subcrit} we will see that the supremum
defining $\mu_q^g$
is attained for $q>\qD$ by a $C^{2,\al}$-spinor. This implies
that the function
  $$[\qD,\infty)\to (0,\infty],\quad q\mapsto \mu_q^g$$
is also continuous from the left.
\end{remark}


\begin{proposition}\label{prop.lamu}
  $$\mu_{\qD}(M,g_0,\si)={1\over\lammin(M,g_0,\si)}$$
\end{proposition}
\proof{}
We have already seen that $\mu_{\qD}$ is conformally invariant, \ie for $g_1\in [g_0]$
  $$\mu_{\qD}(M,g_0,\si)=\mu_{\qD}(M,g_1,\si).$$
The previous proposition states that 
  $$\mu_{\qD}(M,g_1,\si)\geq\mu_2(M,g_1,\si)=\left(\la_1^+(g_1)\right)^{-1}$$ 
if $\vol(M,g_1)=1$, and it follows
  $$\mu_{\qD}(M,g_0,\si)\geq {1\over \lammin(M,g_0,\si)}.$$

It remains to show the inverse inequality which amounts to showing 
  $$\sup_{g\in [g_0]} \la_1^+(g)^{-1}\vol(M,g)^{-1/n}\geq \sup \cF_\qD.$$
For any $\ep>0$ we take a smooth spinor $\phi_\ep$ with 
$\cF_\qD(\phi_\ep)\geq \sup \cF_\qD-\ep$ and $\|\Dir\phi_\ep\|_{L^\qD}=1$.
After a small perturbation of $\phi_\ep$ we can assume that 
$\Dir\phi_\ep$ has no zeros.
We set
  $$g_\ep:=|\Dir\phi_\ep|^{4/(n+1)}g_0.$$
Then we have $\vol(M,g_\ep)=1$, and $|\Dir\phi_\ep|_{g_\ep}$ has constant length $1$. 
Hence,
  $$\cF_\qD^{(M,g_0,\si)}(\phi_\ep)=\cF_\qD^{(M,g_\ep,\si)}(\phi_\ep)
  = \cF_2^{(M,g_\ep,\si)}(\phi_\ep)\leq \mu_2(M,g_\ep,\si)=\la_1^+(g_\ep)^{-1}.$$ 
This implies the proposition.\qed

Now, as we have understood the relationship between the supremum of $\cF_\qD$ and the infimum of $\la_1^+\vol^{1/n}$, 
we want to establish a relationship of the maximizers of $\cF_\qD$ and the minimizers of $\la_1^+\vol^{1/n}$.
This will require some knowledge about the Euler-Lagrange-equation.

One easily sees that 
  $$\cF_q(\alpha \phi+\psi)= \cF_q(\phi)$$
for any $\psi\in \ker\Dir$ and $\alpha\in \mR^*$. 
Hence maximizers of $\cF_q$ appear in families and it will be 
convenient to choose a good representative for each family of maximizers.

\begin{lemma}[Euler-Lagrange equations of $\cF_q$.]\label{lem.euler.lagrange}
Let $q\in [\qD,2]$, and choose $p$ with ${q^{-1}+p^{-1}=1}$. 
Suppose that $\cF_q$ has a maximizing spinor 
$\phi_1\in H_1^q\setminus \ker \Dir$.
Then there is a maximizing spinor $\phi\in H_1^q\setminus \ker \Dir$, $\phi\in \mR^*\phi_1 + \ker D$, 
such that
\begin{equation}\label{eq.nonlin.p}
\Dir\phi= \mu_q^{-1}\, |\phi|^{p-2} \phi,\qquad \phi\in H_1^q, \qquad 
\|\ph\|_{L^p}=1
\end{equation}
\end{lemma}

\proof{} We normalize $\phi_1$ such that $\|\Dir \phi_1\|_{L^q}=1$, thus $\rho_{q,\phi_1}=\mu_q$.
As $\phi_1$ is a maximizer, $d\cF_q(\phi_1):H_1^q\setminus \ker \Dir\to \mR$ is identically zero.
Because of \eref{eq.abl} we see that $\tau:=\phi_1-\mu_q|\Dir\phi_1|^{q-2}\Dir\phi_1\in L^p$ is a weak solution
of $D\tau=0$, hence it is smooth and in the kernel of $\Dir$. Then, $\phi_2:=\phi_1-\tau=\mu_q|\Dir\phi_1|^{q-2}\Dir\phi_1$ satisfies $D\phi_1=D\phi_2$ and hence
\begin{equation}\label{eq.nonlin.q}
  \phi_2 = \mu_q |D\phi_2|^{q-2}D\phi_2, \qquad
  \phi_2\in H_1^q, \qquad\|D\phi_2\|_{L^q}=1.
\end{equation}
Taking norms we obtain 
  $$|\phi_2|=\mu_q|\Dir\phi_2|^{q-1}=\mu_q|\Dir\phi_2|^{q/p}.$$
  $$|D\phi_2|^{q-2}=|D\phi_2|^{-(q-1)(p-2)}=\left(\mu_q^{-1} |\phi_2|\right)^{-(p-2)}$$
Hence $\phi:=\mu_q^{-1}\phi_2$ satisfies \eref{eq.nonlin.p}.
\qed

\begin{theorem}\label{theo.var.corresp}\ \\[-8mm]
\begin{enumerate}[{\rm (a)}]
\item Let $\psi$ be a maximizing spinor of $\cF_\qD$, 
and suppose that $\psi$ is smooth and that $\Dir \psi$ vanishes nowhere. 
Then $g:=|\Dir \psi|^{4/(n+1)}g_0$ is a smooth metric minimizing
$\la_1^+\vol^{1/n}$ in the conformal class $[g_0]$.
\item Let $g\in [g_0]$ be a (smooth) metric minimizing $\la_1^+\vol^{1/n}$,
and let $\psi$ be an eigenspinor of $D_g$ to the eigenvalue $\la_1^+(g)$,
then the length of $|\psi|_g$ is constant and $\psi$ maximizes 
$\cF_\qD$.
\end{enumerate}
\end{theorem}

\proof{}
\begin{enumerate}[(a)]
\item Let $\psi$ be a maximizing spinor. According to the previous lemma, 
there is an $\al\in \mR^*$ and a $\tau\in\ker \Dir$ such that
$\phi:=\al\psi +\tau$ satisfies~\eref{eq.nonlin.p} qith $q=\qD$ and $p=\pD$. 
In particular,
$|\phi|^{\pD-1}=\mu_\qD |\Dir \phi|=\al\mu_\qD|\Dir\psi|$ vanishes nowhere.
The metric $g_\phi:=|\phi|^{4\over n-1}g_0$ satisfies 
$\vol(M,g_\phi)=1$ and $|\phi|_{g_\phi}=1$.  
Then 
  $$\left(\la_1^+(g_\phi)\right)^{-1}\geq \cF_2^{g_\phi}(\phi)=
    \cF_\qD^{g_\phi}(\phi)=
    \cF_\qD^{g_0}(\phi)=
    {\mu_\qD^{-1}\int|\phi|^{\pD}\over \mu_\qD^{-2}
    \left(\int|\phi|^{(\pD-1)\qD}\right)^{2/\qD}}
    =\mu_\qD.$$
As $\mu_\qD^{-1}=\lammin$ we see that $\la_1^+(g_\phi)\leq\lammin$, 
hence $g_\phi$ minimizes $\la_1^+\vol^{1/n}$. By a simple rescaling argument 
one sees that 
$$g:=|\Dir \psi|^{4/(n+1)}g_0=\left({1\over \al}\right)^{4/(n+1)}|\Dir \phi|^{4/(n+1)}g_0=\left({\lammin\over \al}\right)^{4/(n+1)}g_\phi$$
minimizes $\la_1^+\vol^{1/n}$ as well. 
\item By rescaling we can assume that $\vol(M,g)=1$. Unless otherwise
indicated all volume measures, norms, scalar products and 
Dirac operators in this proof are with respect to $g$.
In order to show that $|\psi|_g$ is constant, we define
  $$f_t:={1+t|\psi|_g^2\over \left(\int \left(1+t|\psi|_g^2\right)^n\right)^{1/n}}.$$
One calculates ${d\over dt}|_{t=0}f_t= |\psi|_g^2- \int|\psi|_g^2$.
All metrics $f_t^2 g$ have volume $1$.
Hence, as the infimum of $\la_1^+\vol^{1/n}$ is attained in $g$, we 
have $\la_1^+(g)\leq\la_1^+(g_t)$, and hence
  $$\mu_2(M,g,\si)=\cF_2^g(\psi)\geq \cF_2^{f_t^2g}(\psi)=
    {\int \<D\psi,\psi\>\over \int f_t^{-1}|D\psi|^2}.$$
For $t=0$ equality is attained. Hence 
  $$\int |D\psi|^2\leq \int f_t^{-1} |D\psi|^2.$$
Using $D\psi=\la_1^+(g) \psi$ and deriving with respect to $t$ yields 
  $$0=- \left(\la_1^+(g)\right)^2\int \left(|\psi|_g^2- \int|\psi|_g^2\right)|\psi|_g^2.$$
The right hand side is equal to 
$-\left(\la_1^+(g)\right)^2\left\||\psi|_g^2- \int|\psi|_g^2\right\|_{L^2}=0$,
hence $|\psi|_g$ is constant. 

This implies that $\cF_q^g(\psi)$ is independent of $q$, thus
  $$ \cF_\qD(\psi)=\cF_2^g(\psi)=\left(\la_1^+(g)\right)^{-1}
    =\left(\lammin\right)^{-1}.$$
And hence $\cF_\qD$ attains it supremum in $\psi$. 
\end{enumerate}
\qed

\begin{remark}
Later on, we will see that maximizers of $\cF_\qD$ that vanish nowhere
are always smooth.
\end{remark}

\section{Generalized  metrics}\label{sec.spd}

Unfortunately, we cannot exclude that maximizers of $\cF_q$ vanish 
somewhere. Maximizers with zeros correspond 
to metrics with certain singularities, more precisely to metrics 
whose conformal factor has zeros. 
These metrics are the main object of this section.
We will summarize some facts about the size of the zero set of maximizers 
in Section~\ref{sec.degen}. 

Roughly speaking, generalized metrics are metrics of the form $f^2g_0$ 
where $f\in L^\infty, f\geq 0$. However, for technical and formal 
reasons it is better to use the following definition.

\begin{definition}
Let $g_0$ be a smooth metric on a compact manifold $M$. 
A \emph{generalized metric} is a tuple $(f,g)$ where $g\in [g_0]$ and
$f\in L^\infty$, $f\geq 0$. If $h>0$ is smooth, we identify
$(fh,g)$ with $(f,h^2g)$. Furthermore we identify $g\in[g_0]$ with $(1,g)$.
Generalized metrics having a representative of the form $(f,g_0)$ are called
\emph{conformal} to $g_0$, and those having a representative $(1,g)$, 
$g\in [g_0]$ are called \emph{regular metrics} --- they correspond to metrics 
in the ordinary sense.
The set of all generalized metrics conformal to $g_0$ is called the 
\emph{$L^\infty$-completion} $\ol{[g_0]}$ of the conformal class $[g_0]$.
The \emph{volume of $g\in \ol{[g_0]}$} is defined as $\int f^n \dvol_{g_0}$.
For a generalized metric we say that $f^{-1}(0)$ is the \emph{set of 
degeneration}. 
\end{definition}

\begin{remark}
The reader should pay some attention to the following technical difficulty: 
If $f$ vanishes on an open set, then the 
$L^\infty-(2,0)$-tensor $f^2g_0$ does not determine (the equivalence class) 
$(f,g_0)$.
\end{remark}

However, despite of this remark and slightly abusing the notation, 
we will write $f^2g_0$ instead of $(f,g_0)$. 
Formally $f^2g_0$ is a generalized metric in the above sense, 
not an $L^\infty-(2,0)$-tensor.

Let us now assume that $M$ carries a fixed spin structure $\si$. 
For any generalized metric $g=f^2g_0$ 
we want to define a spinor bundle and a Dirac operator 
on $(M,g,\si)$ in such a way that the results of the previous section 
carry over to this generalization. In particular, the functional $\cF_q$
has to be defined and has to be conformally invariant for $q=\qD$.

As a vector bundle the spinor bundle $\Si(M,g,\si)$ 
is defined to be $\Si(M,g_0,\si)$, and due
to our identification of spinors for different metrics in a fixed conformal
class this construction does only depend on the conformal class $[g_0]$ 
and not on the metric $g_0$ itself. 
For any $x\in M$ and any spinor $\phi$ in the fiber of $\Si(M,g,\si)$ over $x$
we define the pointwise norm
$$|\phi|_g:= \left\{
    \begin{matrix}
     f(x)^{-{n-1\over 2}}|\phi|_{g_0}\hfill & \mbox{if $f(x)\neq 0$}\hfill\cr
     \infty\hfill& \mbox{if $f(x)=0$ and $\phi(x)\neq 0$}\hfill\cr
     0\hfill& \mbox{if $f(x)=0$ and $\phi(x)= 0$}\hfill
    \end{matrix}
\right.
$$
Again, this norm depends only on the conformal class. For smooth sections
$\phi$ of $\Si(M,g,\si)$ such that $|\phi|_g<\infty$ almost everywhere 
we define the $H_1^2(M,g,\si)$-norm as 
  $$\left(\int f^{-1}|\Dir_{g_0}\phi|_{g_0}^2\,\dvol_{g_0}\right)^{1/2}
    + \left(\int |\phi|_{g_0}^\pD\,\dvol_{g_0}\right)^{1/\pD}$$
where we used the conventions $0^{-1}r=\infty$ for $r>0$ and $0^{-1}0=0$. 
The Sobolev space
$H_1^2(M,g,\si)$ is the associated completion. 
\begin{lemma}
There is a natural inclusion 
  $$H_1^2(M,g,\si)\embed H_1^2(M,g_0,\si)$$
\end{lemma}
\proof{}
Cauchy sequences with  respect
to the norm $H_1^2(M,g,\si)$ are also Cauchy sequences with respect to 
$H_1^2(M,g_0,\si)$. Thus we obtain a bounded map  
$H_1^2(M,g,\si)\to H_1^2(M,g_0,\si)$. In order to prove injectivity 
of this map one shows that if $\phi_i$ is a Cauchy sequence in $H_1^2(M,g,\si)$
converging to $0$ with respect to $H_1^2(M,g_0,\si)$, then it converges to $0$
with respect to the $H_1^2(M,g,\si)$ as well.
\qed

We define the Dirac operator 
$\Dir_g:H_1^2(M,g,\si)\to L^2(M,g,\si)$, 
$\Dir_g(\phi):=f^{-1}\Dir_{g_0}(\phi)$.
The spinor $D_g (\phi)$ is well-defined almost everywhere, as
$\{x\in M\,|\,D_{g_0}(\phi)(x)\neq 0 \mbox{ and }f(x)=0\}$ has measure zero.
It is easy to verify that all these definitions only depend 
on the conformal class of $g_0$ and not on $g_0$ itself.
Furthermore, we can reformulate the definition of the above Sobolev space as
\begin{eqnarray*}
 H_1^2(M,g,\si) &=& \Bigl\{\phi \in \Gamma(\Si(M,g,\si))\;
    \Big|\;{\textstyle\int} |\Dir_{g}\phi|_{g}^2\,\dvol_{g}<\infty\Bigr\}.
\end{eqnarray*}

We now extend the definitions $\cF_2^g$ and $\la_1^+(g)$ 
to the $L^\infty$-completion of the conformal class. 
\begin{definition}
For $g=f^2g_0\in \ol{[g_0]}$ and any $\phi\in H_1^2(M,g,\si)\setminus\ker\Dir$ 
we define
  $$\cF_2^g(\phi):={\int \<\Dir_{g_0}\phi,\phi\>_{g_0}\,\dvol_{g_0}
    \over \int f^{-1} |\Dir_{g_0}\phi|_{g_0}^2\,\dvol_{g_0}}.$$
\end{definition}
Because of  
\begin{eqnarray*}
  \int f^{-1} |\Dir_{g_0}\phi|_{g_0}^2\,\dvol_{g_0}
    &\geq& {1\over \|f\|_{L^\infty}}\, \int|\Dir_{g_0}\phi|_{g_0}^2\,\dvol_{g_0}\\
  \cF_2^g(\phi)&\leq & \|f\|_{L^\infty}\,\cF_2^{g_0}(\phi)
\end{eqnarray*}
the functional is well-defined and bounded on $H_1^2(M,g,\si)\setminus \ker D$.

It is also not hard to see that the supremum is attained. In fact, 
let $(\phi_i)$
be a sequence of spinors in  $H_1^2(M,g,\si)\setminus \ker D$ with 
$\cF_2^g(\phi_i)\to \sup \cF_2^g$, 
normalized such that  $\int f^{-1} |\Dir\phi_i|^2=1$. 
Then a subsequence $(\phi_{i_k})$ converges weakly in $H_1^2(M,g,\si)$,
weakly in $H_1^2(M,g_0,\si)$ and strongly in $L^2(M,g_0,\si)$ towards
a $\phi_\infty\in H_1^2(M,g,\si)\setminus \ker D$.
Hence, $\lim \cF_2^g(\phi_{i_k})\leq \cF_2^g(\phi_\infty)$. Thus, 
the supremum is attained in $\phi_\infty$. In analogy to
Lemma~\ref{lem.f.elem} one gets for any smooth 
test spinor $\psi$ 
  $$\int \<f^{-1}\Dir \phi - (\cF_2^g(\phi))^{-1}\,\phi,\Dir \psi\>,$$
which implies 
$\tau:=\Dir_g\phi- (\cF_2^g(\phi))^{-1}\,\phi\in \ker D$ and
finally for $\witi\phi_1:=\phi_\infty+\cF_2^g(\phi)\tau$ 
  $$\Dir_g\witi\phi_1=(\cF_2^g(\witi\phi_1))^{-1}\,\witi\phi_1.$$ 

\begin{definition}
The \emph{first positive Dirac eigenvalue} 
of $(M,g,\si)$, $g=f^2g_0\in \ol{[g_0]}$
is 
  $$\la_1^+(g):=\left(\sup \left\{\cF_2^g(\phi)\,|\,\phi \in H_1^2(M,g,\si)\setminus\ker\Dir\right\}\right)^{-1}.$$
A non-trivial spinor with 
  $$\Dir_g\phi=\la_1^+(g)\phi$$
is an \emph{eigenspinor to the eigenvalue $\la_1^+(g)$}.
\end{definition}

As we have already seen, this definition coincides with the 
definition of the first positive Dirac eigenvalue and a corresponding 
eigenspinor if $g$ is regular. 
  
Most of the statements of the previous section still hold in a 
modified version for generalized metrics, the proofs are nearly identical.
For example one can extend Proposition~\ref{prop.lamu} to
\begin{proposition}\label{prop.bar.nonbar.gen}
  $$\inf_{g\in [g_0])}\la^+_1(g)\vol(M,g)^{1/n}
    =\inf_{g\in \ol{[g_0]}}\la^+_1(g)\vol(M,g)^{1/n} = \mu_{\qD}(M,g,\si)^{-1}.$$
\end{proposition}
The equation $\inf_{g\in [g_0]}\la^+_1(g)\vol(M,g)^{1/n}=\mu_{\qD}(M,g,\si)^{-1}$ 
is exactly the statement of Proposition~\ref{prop.lamu} and the equation 
$\inf_{g\in \ol{[g_0]}}\la^+_1(g)\vol(M,g)^{1/n} = \mu_{\qD}(M,g,\si)^{-1}$ 
can be proven with exactly the same proof.

Similarly, we obtain an analogue of Theorem~\ref{theo.var.corresp}.

\begin{theorem}\label{theo.var.corresp.gen}
Let $(M,g_0,\si)$ be a compact Riemannian spin manifold.\\[-7mm]
\begin{enumerate}[{\rm (a)}]
\item Let $\psi$ be a maximizing spinor of $\cF_\qD$ 
(with regularity $H_1^\infty$).
Then $g:=|\Dir \psi|^{4/(n+1)}g_0$ is a generalized metric minimizing
$\la_1^+\vol^{1/n}$ in $\ol{[g_0]}$, the $L^\infty$-completion of 
the conformal class $[g_0]$.
\item Let $g=f^2g_0\in \ol{[g_0]}$ be a generalized 
metric minimizing $\la_1^+\vol^{1/n}$,
and let $\psi$ be an eigenspinor of $D_g$ to the eigenvalue $\la_1^+(g)$,
then the length of $|\psi|_g$ is constant on $M\setminus f^{-1}(0)$ and 
$\psi$ maximizes $\cF_\qD$.
\end{enumerate}
\end{theorem}

The proof of this theorem is essentially the same as the proof 
of Theorem~\ref{theo.var.corresp}.
Later, we will see that any maximizing spinor has regularity
$C^{2,\al}$, hence $g=|D\psi|^{4/(n+1)}g_0$ is always a generalized metric.

\section{The case of the sphere}\label{sec.sphere}

The sphere $\mS^n=(S^n,\can)$ with the round sphere 
is an important example. On the one hand the invariant $\lammin$
of the round sphere is contained 
in many equations and inequalities. On the other hand the 
round sphere has a large conformal group,
and hence studying the minimizers on the sphere helps to understand 
the analytical difficulties. In particular, we will see
why the conclusion in Theorem~\ref{theo.anydim} does not hold if
$\mu_\qD=\mu_\qD^{\mS^n}$. (Recall $\qD=2n/(n+1)$.)

The invariant $\lammin(\mS^n)$ is not hard to calculate. Recall
that the Yamabe invariant $Y(\mS^n)$ is given by
  $$Y(\mS^n):=\inf_{g\in [\can]}{\int \scal_g\,\dvol_g\over \vol(M,g)^{(n-2)/n}}.$$
It is attained for $g=\can$, and hence $Y(\mS^n)= n (n-1) \, \om_n^{2/n}$,
where $\om_n:=\vol(\mS^n)$.
The Hijazi inequality (\ref{ineq.hij},\ref{ineq.hij.der}) tells us, that
\begin{equation}\label{ineq.sphere.one}
  \lammin(\mS^n)^2\geq {n\over 4 (n-1)}\,Y(\mS^n)= {n^2\over 4}\, \om_n^{2/n}.
\end{equation}
Recall that the sphere of constant sectional curvature $1$ carries
a Killing spinor $\psi$ to the constant $-1/2$, i.e. it satisfies
  $$\na_X \psi= -(1/2) X\cdot \psi.$$
Note that this condition implies 
that the length of $\psi$ is constant. Because of
$D\psi= (n/2) \psi$ we obtain
\begin{equation}\label{ineq.sphere.two}
  \lammin(\mS^n)\leq \la_1^+(\mS^n)\vol(\mS^n)^{1/n}\leq  {n\over 2}\, \om_n^{1/n}.
\end{equation}
And then
  $$\lammin(\mS^n)
    = {n\over 2}\, \om_n^{1/n}.$$
The above Killing spinor satisfies 
  $$\left(\cF_{\qD}(\psi)\right)^{-1}
    = {n\over 2}\, \om_n^{1/n} \left(\sup \cF_{\qD}\right)^{-1},$$
hence $\cF_{\qD}$ attains its supremum in such Killing spinors.
 
Let $A:\mS^n\to \mS^n$ be an orientation preserving 
conformal diffeomorphism. Then the pullback of any spinor 
$\phi$ with respect to $A$ is a section 
of $A^*(\Si(S^n,\can))=\Si(S^n,A^*\can)$, 
and as before we identify $\Si(S^n,A^*\can)\cong \Si(S^n,\can)$ by using 
the map 
$\witi F$ described after Proposition~\ref{prop.conf.change}.
If $\phi$ is a solution of \eref{eq.nonlin.p}, then this pullback, 
denoted by $A^*\phi$, is a solution of this equation as well. 
Furthermore if $\phi$ maximizes $\cF_{\qD}$, then $A^*\phi$ is also a maximizer.
As a consequence, all conformal images
of Killing spinors to the constant $-1/2$ are maximizers of $\cF_{\qD}$.
The following proposition shows that there
are no other maximizers on $\mS^n$.


\begin{proposition}
If $\psi$ is a spinor of regularity $C^2$ that 
attains the supremum of $\cF_{\qD}$, 
then there is a Killing spinor $\phi$ to the Killing constant $-1/2$
and an orientation preserving conformal diffeomorphism $A:\mS^n\to \mS^n$
with $A^*\phi=\psi$.
\end{proposition}

Later on we will see that any maximizer of regularity $H_1^q$ 
with $q>\qD$ is even $C^2$. Hence the statement also holds under this 
weaker assumption.

\begin{lemma}\label{lem.scal.bound}
Let $(M,g,\si)$ be an arbitrary Riemannian spin manifold 
(not necessarily complete or compact). Assume that there is a spinor $\psi$
of constant length $1$ and with $D\psi=\lambda \psi$. Then
  $$\scal = 4  {n-1\over n}\,\lambda^2\,-4 |\witi \nabla \psi|^2,$$
where $\witi\nabla_X \psi:=\na_X\psi + {\la\over n} X\cdot \psi$ denotes
the \emph{Friedrich connection} on spinors.
\end{lemma}

\proof{of the lemma}
The Friedrich connection is a metric connection, hence
  $$0={1\over 2}\,d^*d\<\psi,\psi\> = {\rm Re} \<\witi\na^* \witi\na \psi,\psi\> 
  - \<\witi\na\psi,\witi\na\psi\>.$$
The twisted version of the Schr\"odinger-Lichnerowicz formula 
yields
  $$\left(D-{\lambda\over n}\right)^2=\witi\na^*\witi\na + {\scal\over 4}- {(n-1) \lambda^2 \over n^2}.$$
Hence 
  $$\left(\la-{\la\over n}\right)^2= \<\witi\na^*\witi\na \psi,\psi\>+ {\scal\over 4}- {(n-1) \lambda^2 \over n^2}.$$
We obtain
  $$ {n-1\over n} \lambda^2
   = \<\witi\na\psi,\witi\na\psi\> + {\scal\over 4}.$$
\qed
\proof{of the proposition}
 

We have to show that the supremum is not attained by any other spinor. 
For proving this, assume that $\psi$ is a maximizer, 
$\|\psi\|_{L^\qD}=1$. 
As the Dirac operator on $\mS^n$ has kernel $\{0\}$, the spinor 
satisfies the Euler-Lagrange 
equation \eref{eq.nonlin.p}.

On the open subset $S^n\setminus \psi^{-1}(0)$ we define the metric 
  $$g_1:= |\psi|^{4\over n-1}\can.$$
In this metric \eref{eq.nonlin.p} transforms into a solution of 
  $$D_{g_1}\psi=\lammin \psi\qquad |\psi|_{g_1}\equiv 1.$$

On the one hand, one calculates 
\begin{eqnarray*}
   \scal_{g_1}&=& 
    4 {n-1\over n-2}\,|\psi|^{-{n+2\over n-1}}\Delta_{\can} |\psi|^{n-2\over n-1}
    + \scal_{\can}|\psi|^{-{4\over n-1}}\\
\end{eqnarray*}
and integration yields
\begin{eqnarray}
\int_{S^n\setminus \psi^{-1}(0)}\scal_{g_1}\,\dvol_{g_1}
    &=& \int_{M\setminus \psi^{-1}(0)} 4 {n-1\over n-2}\,|\psi|^{{n-2\over n-1}}\Delta_{\can} |\psi|^{n-2\over n-1}
    + n(n-1)|\psi|^{2{n-2\over n-1}}\,\dvol_{\can}\nonumber\\
&=& \int_{S^n}
    4 {n-1\over n-2}\,\left|d|\psi|^{n-2\over n-1}\right|_{\can}^2
    + n(n-1)|\psi|^{2{n-2\over n-1}}\,\dvol_{\can},\label{ineq.preyamb}
\end{eqnarray}
where the last term arises by partial integration. In order to make this step 
precise one has to exhaust $M\setminus \psi^{-1}(0)$ by smooth manifolds
with boundary and partially integrate over these exhausting manifolds. 
The boundary terms vanish in the limit, as $\psi\to 0$ on the boundaries and $d|\psi|$ is bounded.

It is a standard fact from the resolution of the Yamabe problem (see e.g.\
\cite{lee.parker:87})
that any $H_1^2$ function $f$ satisfies
\begin{eqnarray}\label{ineq.yamafunc}
{\int_{S^n}
    4 {n-1\over n-2}\,\left|df\right|^2
    + n(n-1)f^2\,\dvol_{\can}\over \left(\int_{S^n} f^{2n\over n-2}\,
    \dvol_{\can}\right)^{n-2\over n}} &\geq& Y(\mS^n)=n(n-1)\om_n^{2/n}.
\end{eqnarray}
Setting $f:=|\psi|^{n-2\over n-1}$, we obtain $\int_{S^n} f^{2n\over n-2}=1$, 
and hence the right side of \eref{ineq.preyamb} is bounded from below 
by $n(n-1)\om_n^{2/n}$.

One the other hand, the previous lemma provides
  $$\scal_{g_1}\leq 4 {n-1\over n}\left(\lammin\right)^2=n(n-1)\om_n^{2/n},$$
and as  $\vol(S^n\setminus \psi^{-1}(0),g_1)=1$, we see that we must have
equality in all inequalities involved, in particular in~\eref{ineq.yamafunc}. 
An application of the maximum principle 
\cite{lee.parker:87} yields that $f$ does not vanish, and hence $\psi$ 
has no zeros. Furthermore, $g_1$ is a metric of constant scalar curvature
conformal to $\mS^n$, and such a metric is necessarily of the form
$g_1:=A^*\can$ for an orientation preserving conformal diffeomorphism
$A:\mS^n\to \mS^n$. With respect to $g_1$ we obtain 
$\witi \na \psi=0$, hence $\psi$ is a Killing spinor on $(S^n,g_1)$. 
This implies that $\phi:=(A^{-1})^*\psi$ is a Killing spinor on $\mS^n$.
\qed


\begin{remark}
There are solutions to \eref{eq.nonlin.p} that do not maximize the functional.
An easy construction of such a solution is as follows. 
The map $A:\mS^2=\mC\cup\{\infty\}\to \mS^2=\mC\cup\{\infty\}$, 
$z\mapsto z^k$, $k\in\mN\setminus\{0\}$ 
is conformal with branching points $0$ and $\infty$. If $\phi$ is a solution
of $D\phi=c|\phi|^2\phi$, then 
the pullback $\psi:=A^*\phi$  is a solution  of 
$D\psi=c|\psi|^2\psi$ on $\mS^2\setminus\{0,\infty\}$.
Here $\psi$ is a section of the pull-backed spinor bundle, which is defined 
using the pull-backed spin structure. The pull-backed spin structure
 on $\mS^2\setminus\{0,\infty\}$ coincides with the standard spin structure
iff $k$ is odd. If $k$ is odd, one can show that the extension of $\psi$ by 
setting $\psi(0)=0$ and $\psi(\infty)=0$ is a solution of 
$D\psi=c|\psi|^2\psi$ on $\mS^2$.
However, $\int|\psi|^4= k \int|\phi|^4$. This implies 
$\cF(\psi)=k^{-1/3}\cF(\phi)$. Hence, if $\phi$ is a Killing spinor, 
and $k\geq 3$, $k$ odd, then $\psi$ is a non-maximizing solution.
\end{remark}

\section{Regularity theorems for the Euler-Lagrange equations}\label{sec.reg.theo}

In this section we want to study solutions of the Euler-Lagrange equations.
The first subsection ``Removal of singularities'' will be used in 
the following section to 
extend a solution of $\mR^n=\mS^n\setminus\{\mbox{South Pole}\}$ 
to a solution on $\mS^n$. The second subsection states that solutions
of the Euler-Lagrange equations are $C^{2,p-2}$.

\subsection{Removal of singularities}

\begin{theorem}[Removal of singularities theorem]\label{theo.sing.rem}
Let $p\in[{n\over n-1},\infty)$.
Let $(U,g)$ be a (not necessarily complete) Riemannian manifold 
equipped with a spin structure, let $x\in U$.
Assume that $\phi\in L^p(\Si(U\setminus\{x\}),g)$ 
satisfies weakly on $U\setminus\{x\}$ the equation
\begin{equation}\label{eq.singsol}
    D\phi = \la |\phi|^{p-2}\phi.
\end{equation}
Then this equation even holds weakly on $U$. In particular, the distribution
$D\phi$ does not have singular support in $x$ and is contained in 
$L^q$. 
\end{theorem}

\proof{}
Let $\psi$ be a smooth spinor compactly supported in $U$.
We have to show 
\begin{eqnarray}\label{eq.singsol.show}
  \int_U\<\phi,\Dir \psi\>= \la \int_U \< |\phi|^{p-2}\phi,\psi\>.
\end{eqnarray}

For any small $\ep>0$ we choose a smooth 
cut-off function $\eta_\ep:U\to[0,1]$
with $\eta_\ep\equiv 1$ on $B_\ep(x)$, with $|\na \eta_\ep|\leq 2/\ep$
and with support in $B_{2\ep}(x)$.
We rewrite the left hand side as
\begin{eqnarray}
         \int_U\<\phi,\Dir \psi\>
         &=&\int_U\Bigl\<\phi,\Dir \Bigl( (1-\eta_\ep) \psi + \eta_\ep \psi\Bigr)\Bigr\>\\
         &=&\int_U\Bigl\<\phi,\Dir \Bigl( (1-\eta_\ep)\Bigr) \psi\Bigr\> 
            + \int_U\<\phi,\eta_\ep \Dir \psi\> + \int_U\<\phi,\na \eta_\ep \cdot \psi\>\nonumber
\end{eqnarray} 
As $\phi$ is a weak solution of \eref{eq.singsol} 
on $U\setminus\{x\}$, the first term
equals to 
  $$\la \int \< |\phi|^{p-2}\phi,(1-\eta_\ep)\psi\>,$$ 
and for $\ep\to 0$ it tends to the right hand side of~\eref{eq.singsol.show}.

Let $q$ be related to $p$ via $1/q+1/p=1$.
The absolute value of the second term is bounded by 
  $$\|\phi\|_{L^p(B_{2\ep})}\,\|D\psi\|_{L^q(B_{2\ep})}$$
which tends to $0$ for $\ep\to 0$.

Finally, the absolute value of the third term is bounded by
  $${2\over \ep}\,\|\phi\|_{L^p(B_{2\ep})}\,\|\psi\|_{L^q(B_{2\ep})}\leq C\,\|\phi\|_{L^p(B_{2\ep})}\,
    \ep^{{n\over q}-1}.$$
Our condition $p\geq {n\over n-1}$ yields $q\leq n$, and hence the third term 
also tends to $0$ for $\ep\to 0$.
\qed

\subsection{Regularity}

\begin{theorem}[$C^{1,\alpha}$-regularity theorem]\label{theo.reg}
Suppose that $\phi\in H_1^q$, $q\in [\qD,2]$
is a solution of equations~\eref{eq.nonlin.p}.
Suppose that there is an $r>\pD$ such that $\|\phi\|_{L^r}<\infty$.
Then $\phi$ is $C^{1,\alpha}$ for any $\alpha\in (0,1)$.\\
Furthermore, we obtain a uniform bound of the $C^{1,\alpha}$-norm in the 
following sense.
Let us choose $k,K>0$ such that $\|\phi\|_{L^r}<k$ and $\mu_q\geq K$.
Then for any $\al\in (0,1)$ there is a constant
 $C$ depending only on $(M,g,\si)$, $p$, $r$, $K$, $k$ and $\al$ with
  $$\|\phi\|_{C^{1,\al}}\leq C.$$
\end{theorem}

\begin{remark}The following
example shows that the theorem cannot hold without the $L^r$-bound. 
Let $M=\mS^n$ and $p=\pD$. 
Let $\psi$ be a Killing spinor to the Killing constant $-1/2$.
Suppose $\|\psi\|_{L^\pD}=1$. 
Let $A:\mS^n\to \mS^n$ be a M\"obius transformation, 
such that the differential in $x\in \mS^n$ 
satisfies $(dA)_x=2 \Id$.
Then $\psi_i=d(A^i)\psi$ is a family of solutions of \eref{eq.nonlin.p}, 
maximizing $\cF_\qD$.
However, for any $r>\pD$ one can show that $\|\psi_i\|_{L^r}\to \infty$ 
for $i\to \infty$.
Hence, the $L^r$-bound is necessary for the theorem to hold. 
\end{remark}

\begin{remark}
The theorem will be applied in several versions. At first, we will apply it when
$p<\pD$. In this case the Sobolev embedding already provides the required 
$L^r$-bound on $\phi$. However, the $r$ given by the Sobolev embedding depends on $p$. It will
be of central importance to obtain a bound that is uniform for $p\to \pD$. 
After having proved Theorem~\ref{theo.anydim}  
the uniformity statement in the above theorem will be used to 
obtain a bound that is uniform for $p\to \pD$.
Finally, the regularity theorem will be applied in the case $p=\pD$. In this case an additional $L^r$-bound 
is required as well. 
\end{remark}

\proof{}
The proof uses the following ``bootstrap argument''.
At first, we assume $r<n\,{n+1\over n-1}$.
As $\phi$ is $L^r$, the right hand side of  \eref{eq.nonlin.p}, i.e.\ $|\phi|^{p-2}\phi$, 
is $L^{r/(p-1)}\embed L^s$ with $s:={r/(\pD-1)}= r {n-1\over n+1}< n$.
We apply the Global $L^p$-estimates~\ref{theo.global.lp} and get $\phi\in H_1^s$.
Using the Sobolev embedding I, Theorem~\ref{theo.sobo} (a)
one obtains $\phi\in L^{r'}$ with $r'={ns\over n-s}={rn {n-1\over n+1}\over n-r{n-1\over n+1}}$.
Using $r>\pD=2n/n-1$ one sees that $r'>r$, hence we have obtained stronger regularity for $\phi$.
We iterate this argument and get $L^{\ti r}$-bounds for arbitrarily 
large~${\ti r}$. For any
$\ti r>n\,{n+1\over n-1}$, we obtain $\phi\in H_1^{\ti s}$ with $\ti s:={\ti r/(\pD-1)}>n$. We apply 
the Sobolev embedding theorem~II, Theorem \ref{theo.sobo} (c) and obtain $\phi\in C^{0,\al}$ for any $\al\in (0,1)$.
Hence $|\phi|^{p-2}\phi$ is $C^{0,\al}$ as well, and applying  Schauder estimates \ref{theo.schauder}
we get $\phi\in C^{1,\al}$ for arbitrary $\al$.

The uniformity of the upper bound is clear from the construction.
\qed

The bootstrap can be continued and we obtain better regularity.

\begin{proposition}[Improved regularity]\label{prop.imp}
We assume the assumptions of the previous theorem.
\begin{enumerate}[{\rm (1)}]
\item Let $U:=M\setminus\phi^{-1}(0)$. Then $\phi|_U\in C^{\infty}(U)$.
\item If $p > 2$, 
then $\phi\in C^{2,\alpha}$  for any $\alpha\in (0,1)\cap(0, p-2]$.
Furthermore,
  $$\|\phi\|_{C^{2,\al}}\leq C,$$
where $C$ depends only on $(M,g,\si)$, $p$, $r$, $K$, $k$ and $\al$.
\item If $n=2$ and $p=\pD=4$, then $\phi\in C^\infty$.  
Furthermore,
  $$\|\phi\|_{C^m}\leq C,$$
where $C$ depends only on $(M,g,\si)$, $p$, $r$, $K$, $k$ and $m$.
\end{enumerate}
\end{proposition}

\proof{} (1) On $U$  we can continue the
bootstrap argument and apply inductively the Schauder estimates
Theorem~\ref{theo.schauder}.
We conclude that $\phi$ is smooth on $U$. We obtain~(1). 
However, uniform bounds will be difficult to obtain.

(2) The case $p=2$ is trivial. 
Let $p>2$. We know that $\phi$ is $C^{1,\al}$ for any $\alpha$.
Hence, using Appendix~\ref{sec.scha.power} one sees that 
$|\phi|^{p-2}\phi$ is $C^{1,\al}$. The Schauder estimates imply that $\phi$
is $C^{2,\al}$.

(3) Similarly, if $n=2$ and $p=4$, 
then $\phi\mapsto |\phi|^{p-2}\phi$ is also smooth in $0$. Hence the bootstrap 
can go on with higher order Schauder
estimates, and we inductively get $\phi\in C^m(M)$ for any $m$. 
The construction obviously provides uniform bounds. We obtain~(3).
\qed

\section{Solution of the variational principle}\label{sec.solution}

\begin{proposition}\label{prop.sol.subcrit}
Let $q$ and $p$ be related via $1/p+1/q=1$. For $q\in(\qD,2)$ 
the supremum $\mu_q$ of $\cF_q$ is attained by a spinor field
$\phi\in C^{2,\al}$, $\alpha\in (0,1)\cap(0, p-2]$, 
$1/p+1/q=1$. The spinor $\phi$ can be chosen such that $\phi$
is a solution of \eref{eq.nonlin.p}.
\end{proposition}

\proof{}
Let $\phi_i$ be a maximizing sequence for $\cF_q$, \ie
$\cF_q(\phi_i)\to \mu_q$. We may assume $\|\Dir \phi_i\|_{L^q}=1$, and that 
$\phi_i$ is orthogonal to $\ker \Dir$.
After taking a subsequence there is a $\phi_\infty\in H_1^q$ 
such that $\phi_i$ converges weakly to $\phi_\infty$ in $H_1^q$.
Thus $\|\Dir\phi_\infty\|_{L^q}\leq \liminf \|\Dir\phi_i\|_{L^q}=1$.
The compactness of the embedding $H_1^q\hookrightarrow L^p$ provides
a subsequence, that converges strongly to
$\phi_\infty$ in $L^p$. 
This implies
  $$\int\<D\phi_i,\phi_i\>=
    \underbrace{\int\<D\phi_i,\phi_i-\phi_\infty\>}_{\leq\|D\phi_i\|_{L^q}\|\phi_i-\phi_\infty\|_{L^p}} 
    +\underbrace{\int\<D\phi_i,\phi_\infty\>}_{\to\int\<D\phi_\infty,\phi_\infty\>}\to \int\<D\phi_\infty,\phi_\infty\>.
$$
Hence,
  $$\mu_q= \lim \cF_q(\phi_i)= \lim
    {\int\<\Dir \phi_i,\phi_i\>\over \|\Dir\phi_i\|_{L^q}^2}\leq 
    {\int\<\Dir \phi_\infty,\phi_\infty\>\over \|\Dir\phi_\infty\|_{L^q}^2}=
    \cF_q(\phi_\infty)\leq \mu_q.$$
As a consequence, we have equality in all inequalities, in particular
$\|\Dir\phi_\infty\|_{L^q}=1$.
According to Lemma~\ref{lem.euler.lagrange} one can 
find $\al\in \mR^*$ and $\tau\in \ker \Dir$ such that 
$\phi=\al\phi_\infty+\tau$
solves \eref{eq.nonlin.p}.
Proposition~\ref{prop.imp} (Improved Regularity) 
tells us that $\phi$ is $C^{2,\al}$.
\qed

\begin{theorem}[Uniform $C^0$-estimate]\label{theo.anydim}
Let $\phi$ be a solution of \eref{eq.nonlin.p} with $p\in [2,\pD)$
and $\mu_q \geq \mu_{\qD}^{\mS^n}+\ep$, $\ep>0$.
Then there is a constant $C=C(M,g,\si,\ep)$ such that
  $$\|\phi\|_{C^0} < C.$$
\end{theorem}

\begin{remark}\label{rem.sphereconc}
The conclusion of this theorem does not hold any longer, if we drop 
the condition $\mu_q>\mu_\qD^{\mS^n}$. In fact, if $q=\qD$ and 
$(M,g,\si)=\mS^n$, then there is the following counterexample: 
Let $A$ be a non-isometric, but conformal map $\mS^n\to\mS^n$ 
fixing North and South pole. Let $\phi$ be a Killing spinor. As seen in the 
previous section $\phi$ is a maximizer of $\cF_\qD$ and a solution
of \eref{eq.nonlin.p}. The images $\phi_k:=(A^k)_*\phi$ under the conformal 
maps $A^k$ are also  maximizers of $\cF_\qD$ and solutions
of \eref{eq.nonlin.p}. However, as easily seen, $\|\phi_k\|_{C^0}$
tends to $\infty$ for $k\to \infty$.
\end{remark}

\proof{of the theorem}
Assume that such a constant does not exist. Then we find a sequence
of solutions $\phi_k$ of
\begin{equation}\label{eq.nonlin.pk}
    \Dir\phi_k= \mu_{q_k}^{-1}\, |\phi_k|^{p_k-2} \phi_k,\qquad \phi_k\in H_1^{q_k}, \qquad 
\|\phi_k\|_{L^{p_k}}=1
\end{equation}
$1/q_k = 1 - 1/p_k$, $\mu_{q_k} \geq \mu_{\qD}^{\mS^n}+\ep$ and
\begin{eqnarray}\label{form.div}
\|\phi_k\|_{C^0}\to \infty.
\end{eqnarray}

Let us assume for a moment that $p_\infty:=\liminf p_k < \pD$.
In this case, we can choose a subsequence with $p_k\to p_\infty$.
We have 
$1=\|\phi_k\|_{L^{p_k}}^{p_k}= \Bigl\|\,|\phi_k|^{p_k-2}\phi_k\Bigr\|_{L^{q_k}}^{q_k}={\mu_{q_k}}^{q_k}\|D\phi_k\|_{L^{q_k}}^{q_k}$. 
We conclude that 
$\phi_k$ is bounded in $H_1^{\ti q}$ for a $\ti q > \qD$,
and hence in $L^r$ for an $r>\pD$.
Then the regularity theorem (Theorem~\ref{theo.reg}) says that
$\|\phi_k\|_{C^0}$ is bounded, in contradiction to \eref{form.div}.
Hence, the case $p_\infty>\pD$ can not occur, i.e.\ $\lim p_k = \pD$.

There is a sequence of points $s_k\in M$ with
  $$m_k:=|\phi_k(s_k)|=\max\left\{|\phi_k(x)|\,\big|\,x\in M \right\}\to \infty.$$
The idea is to blow up suitably the metric such that we obtain
in the limit a solution on Euclidean $\mR^n$. 

We define
\begin{eqnarray*} 
  \ti g_k&:=&(m_k)^{2(p_k-2)}g\\
  \ti \phi_k&:=& (m_k)^{(p_k-2)\frac{n-1}{2}-1} \phi_k.
\end{eqnarray*}
One easily verifies
\begin{equation}\label{eq.normed}
   |\ti \phi_k(s_k)|_{\ti g_k}=1,
\end{equation}
and obviously \eref{eq.nonlin.p} transforms into
\begin{equation}\label{eq.preser}
   D_{\ti g_k}\ti \phi_k= \frac{1}{\mu_{q_k}}|\ti\phi_k|_{\ti g_k}^{p_k-2}\ti\phi_k. 
\end{equation}
We calculate
\begin{eqnarray*} 
\left\|\ti\phi_k\right\|_{L^{p_k}(\ti g_k)}^{p_k} 
   &= &\mu_{q_k}\int\<D_{\ti g_k}\ti\phi_k,\ti\phi_k\>_{\ti g_k}\,\dvol_{\ti g_k}\\
   &=& \mu_{q_k}\int\<D_g\ti\phi_k,\ti\phi_k\>\,\dvol_g\\
   &=&\left(m_k\right)^{2 \left((p_k-2)\frac{n-1}{2}-1\right)}
   \underbrace{\left\|\phi_k\right\|_{L^{p_k}(g)}^{p_k}}_{=1}.  
\end{eqnarray*}
We can assume that $m_k\geq 1$. 
As $p_k\leq\pD$ implies $(p_k-2)\frac{n-1}{2}-1\leq 0$, we obtain
\begin{equation}\label{ineq.unif.lp}
  \left\|\ti\phi_k\right\|_{L^{p_k}(\ti g_k)}\leq 1.
\end{equation}
The injectivity radius of $(M,\ti g_k)$ tends to infinity, i.e.\ 
for any $R>0$ there is a $k_0=k_0(R)\in\mN$ 
such that for all $k\geq k_0$ the exponential map 
$\exp^{\ti g_k}_{s_k}:T_{s_k}M\to M$ 
with respect to $\ti g_k$ and based in $s_k$ is a diffeomorphism 
on the ball of radius $R$ around $0$.
Now, we identify $(T_{s_k}M,\ti g_k)$ with $(\mR^n,\geucl)$.
Then $\overline g_k:=(\exp^{\ti g_k}_{s_k})^*{\ti g_k}$ 
is Riemannian metric on $B_R(0)$ that
coincides with $\geucl$ in $0$.
In the limit $k\to \infty$ the metrics $\overline g_k$ converge 
to $\geucl$ in the $C^\infty$-topology.

As already said in previous sections, 
the construction of the spinor bundle, its scalar product 
and its connection depends on the Riemannian metric (and the spin structure).
In order to work out the blowup construction one has to define
a pull-back of the spinors 
$\witi\phi_k\in \Gamma(\Si(M,g_k,\si))$ via normal coordinates and then obtain
a spinor in $\Gamma(\Si(B_R(0),\geucl))$. 
In the literature two such pull-backs are used, a pullback
construction carried out in \cite{ammann.humbert.morel:p03av2} and inspired by
\cite{bourguignon.gauduchon:92}, or via radial parallel transport 
(see e.g.\ the solution of the index problem in \cite{roe:88} 
using Getzler rescaling). Both pullbacks can be used here.
The pullback in \cite{ammann.humbert.morel:p03av2} has better 
approximation properties, and is an important tool if one wants to study
fine asymptotics of the blowup. However, the pullback via 
radial parallel transport is technically simpler to introduce, 
hence it will be used here. 

For $R>0$ let $\Si_0(B_R(0),\geucl)$ (resp. $\Si_{s_k}(M,\overline g_k,\sigma)$)
be the fiber of $\Si(B_R(0),\geucl)$ over $0$ 
(resp.  $\Si(M,\overline g_k,\sigma)$ over $s_k$).
Let  us define the radial vector
field $X=r{\partial r}=\sum x^i\partial{x^i}$ on $B_R(0)\subset \mR^n$. Its 
length is the distance from $0$.
For sufficiently large $k$, the exponential map of $(M,\overline g_k)$ based 
in $s_k$, denoted by $\exp^{\ti g_k}_{s_k}$ is a diffeomorphism from
$B_R(0)$ onto its image.
One chooses a (complex) linear isometry 
$\Si_0(B_R(0),\geucl)\to\Si_{s_k}(M,\overline g_k,\sigma)$.
This map extends uniquely to a fiber preserving map 
$A:\Si(B_R(0),\geucl)\to\Si(M,\overline g_k,\sigma)$, such that
  $$\begin{matrix}
    \Si(B_R(0),\geucl)& \stackrel{A}{\longrightarrow} & \Si(M,\overline g_k,\sigma)\cr
   \downarrow & & \downarrow\cr
   B_R(0) &  \stackrel{\exp^{\ti g_k}_{s_k}}{\longrightarrow} & M
 \end{matrix}
 $$
commutes and such that
  $$A(\nabla_X \phi)= \nabla_{(\exp^{\ti g_k}_{s_k})_*(X)}A\phi.$$
In a neighborhood of $s_k$, this condition can be equivalently 
characterized by saying that 
$\phi\mapsto A\circ \phi \circ (\exp^{\ti g_k}_{s_k})^{-1}$ maps
parallel sections of $\Si(B_R(0),\geucl)$ to radially parallel 
sections of $\Si(B_R(s_k),\overline g_k,\si)$.
One easily sees that $A$ is a fiberwise isometry, but the connection 
on the spinor bundle is not preserved. However, this will not matter, as 
in the limit $k\to \infty$, $R$ fixed, the connections converge in the 
$C^\infty$-topology.

Let $k\geq k_0(R)$. On $(B_R(0),\geucl)$ we define the spinor  
  $$\overline \phi_k:=A^{-1}\circ \witi\phi_k\circ \exp^{\ti g_k}_{s_k}$$
and the operator $D_k:\Gamma(\Si(B_R(0),\geucl))\to\Gamma(\Si(B_R(0),\geucl))$,
  $$D_k:=A^{-1}\circ D_{\overline g_k}\circ A,$$
where 
$D_{\ol g_k}$ is the Dirac operator on $(B_R(0),\overline g_k)$.
We obtain
  $$D_k\overline{\phi}_k= {1\over \mu_{q_k}}\,|\overline{\phi}_k|^{p_k-2}\overline{\phi}_k$$
Note that 
  $$\|\overline\phi_k\|_{C^0(B_R(0))}\leq |\overline\phi_k(0)|=1.$$
Hence, we may apply the interior $L^p$- and the Schauder-estimates 
\ref{theo.int.lp} and \ref{theo.schauder}
to conclude that
  $$\|\overline{\phi}_k\|_{C^{1,\al}(B_{R/2}(0))}\leq C(R),$$
with constants $C(R)$ and $k(R)$.
The constant $C(R)$ does not depend on $k$ if $\overline g_k$ is sufficiently
$C^\infty$-close to $\geucl$. In particular, the $C^\infty$-convergence
$\overline g_k\to \geucl$ says that $C(R)$ does not depend on 
$k$ for $k\geq k_1(R)$.

Compare $D_k$ with the Dirac operator $\Dflat$ on Euclidean $\mR^n$.
(For similar and more explicit calculations the reader might consider 
e.g.\ \cite{pfaefflediss} for pullback via radial parallel transport,
and  \cite{ammann.humbert.morel:p03av2} for the other pullback method). 
We have
  $$\|(\Dflat -D_k)\overline\phi_k\|_{C^{0,\al}(B_{R/2}(0),\geucl)}\leq \tau_{k,R} \|\overline\phi_k\|_{C^{1,\al}(B_{R/2}(0),\overline g_k)},$$
with $\lim_{k\to \infty}\tau_{k,R}=0$. 
The convergence is not uniform in $R$,
but this will not matter in the following.

We choose a sequence of radii $R_m\to \infty$. 
For each $R_m$, the Arcela-Ascoli theorem (Theorem~\ref{theo.arc.asc}) 
allows us to choose a 
subsequence of $(\overline\phi_k)$ converging in 
$C^1(B_{R_m}(0),\geucl)$.
After passing to a diagonal sequence, we see that there is a spinor
$\overline{\phi}_\infty$ on $\mR^n$, such that 
$\res{\overline{\phi}_k}{B_{R}(0)}$
converges to 
$\res{\overline{\phi}_\infty}{B_{R}(0)}\in C^1(B_R(0),\geucl)$ 
for all $R>0$.

Then $\overline\phi_\infty$ is a solution of
  $$\Dflat \overline\phi_\infty = {1\over\mu_\qD}\, |\overline\phi_\infty|^{\pD-2}
     \overline\phi_\infty$$
on $\mR^n$.

The estimate~\eref{ineq.unif.lp} says that
$\|\witi\phi_k\|_{L^{p_k}(B_R(s_k),\witi g_k)}\leq 1$. Hence,
for any $\ep>0$ and $R>0$ there is $k_2=k_2(R,\ep)$ such that
  $$\|\overline\phi_k\|_{L^{p_k}(B_R(0))}\leq 1 + \ep$$
for all $k\geq k_2$.
Because of the $C^1$-convergence $\overline\phi_k\to \overline\phi_\infty$,
Fatou's lemma yields
  $$\|\overline\phi_\infty\|_{L^{\pD}(B_R(0))}\leq 1$$
for any $R\in (0,\infty)$, and finally for $R=\infty$.

We identify $\overline\phi_\infty$ via stereographic projection with an
$L^\pD$-spinor $\widehat\phi_\infty$ on $\mS^n\setminus\{\mbox{South pole}\}$
with the identification provided by the application $\witi F$ directly after 
Proposition~\ref{prop.conf.change}. We obtain
\begin{eqnarray}\label{eq.Sn}
  D^{\mS^n} \widehat\phi_\infty = {1\over \mu_\qD}\, |\widehat\phi_\infty|^{\pD-2}
     \widehat\phi_\infty
\end{eqnarray}
and $\|\widehat \phi_\infty\|_{L^\pD}\leq 1$. 
The removal of singularities theorem, i.e.\ Theorem~\ref{theo.sing.rem}, 
says that \eref{eq.Sn} holds on the whole sphere $\mS^n$.

\begin{eqnarray*}
   \int_{\mS^n}\<D\widehat\phi_\infty,\widehat\phi_\infty\> &=&
         \mu_\qD^{-1}\,\|\widehat\phi_\infty\|_{L^\pD(\mS^n)}^\pD,\\
  \|D\widehat\phi_\infty\|_{L^\qD(\mS^n)}&=& 
         \mu_\qD^{-1}\,\Bigl\|\,|\widehat\phi_\infty|^{\pD-1}\widehat\phi_\infty\Bigr\|_{L^\qD(\mS^n)}
        =   \mu_\qD^{-1}\,\|\widehat\phi_\infty\|_{L^\pD(\mS^n)}^{\pD-1},\\
   \mu_\qD^{\mS^n}\geq \cF^{\mS^n}_\qD(\widehat\phi_\infty)&=& 
         \mu_\qD\|\widehat\phi_\infty\|_{L^\pD(\mS^n)}^{2-\pD}\geq \mu_\qD
\end{eqnarray*}
which is apparently a contradiction to our assumption
$\mu_\qD\geq\mu_\qD^{\mS^n}+\ep$.
\qed

\begin{proposition}
If there is a $p_0<\pD$ and an $r>\pD$ such that for all
$t\in (p_0,\pD)$ there is a solution $\phi_t$ of equation~\eref{eq.nonlin.p}
with $p=t$, $1/q+1/p=1$ and such that $\|\phi_t\|_{L^{r}}$ is bounded
by a constant $C$ independent from $t$,
then there is a sequence $t_i\to \pD$ such that $\phi_{t_i}$ converges
in the $C^1$-topology
to a solution of equation~\eref{eq.nonlin.p} with $p=\pD$.
\end{proposition}

\proof{}
For $p$ sufficiently close to $\pD$, we know because of
Proposition~\ref{prop.func}
that $\mu_t$ is bounded from below by a positive constant. Thus,
we can apply the regularity theorem (Theorem~\ref{theo.reg}) which tells us
that $(\phi_t)$ is uniformly bounded
in $C^{1,\al}$. 
Hence, for a  sequence $(t_i)$ with $t_i<\pD$,
converging to $\pD$,
the spinor fields $\phi_{t_i}$ converge in the $C^1$-topology to a
$C^1$-spinor field $\phi_{\pD}$ which is a solution of 
equation~\eref{eq.nonlin.p} with $p=\pD$. \qed

\section{Proof of the main theorem}\label{sec.proofmain}

In this section we want to prove the main result of this publication, namely
Theorem~\ref{theo.main}.

\proof{of Theorem~\ref{theo.main}}
Proposition~\ref{prop.sol.subcrit} tells us that
for any $q\in (\qD,2)$, $\qD:=2n/(n+1)$, the functional $\cF_q$ is attained by a maximizer
denoted $\phi_q$ satisfying 
$$
\Dir\phi_q= \mu_q^{-1}\, |\phi_q|^{p-2} \phi_q,\qquad \phi_q\in H_1^q, \qquad 
\|\phi_q\|_{L^p}=1
$$
where $p$ and $q$ are related via $p^{-1}+q^{-1}=1$.
We have assumed that $\lammin(M,g,\si)< \lammin(\mS^n)$ 
which is equivalent to $\mu_\qD>\mu_\qD^{\mS^n}$.
As the function $q\mapsto \mu_q$ is continuous from the right
(Proposition~\ref{prop.func}), we see that there is an $\ep>0$ such that
 $\mu_q>\mu_\qD^{\mS^n}+\ep$ for $q$ close to $\qD$.
For such $q$, Theorem~\ref{theo.anydim} implies that $\phi_q$ are uniformly bounded in the $C^0$-norm, 
and then we can use Theorem~\ref{theo.reg} to conclude that these $\phi_q$ are even uniformly bounded 
in $C^{1,\al}$. The Theorem of Arcela-Ascoli 
(Theorem~\ref{theo.arc.asc}) 
implies that there is there is a sequence $q_i\to \qD$ such 
that $\phi_{q_i}$ converges in the $C^1$-norm 
to a solution $\phi$ of 
$$
\Dir\phi= \lammin\, |\phi|^{\pD-2} \phi,\qquad \phi\in C^1, \qquad 
\|\phi\|_{L^\pD}=1
$$
Theorem~\ref{theo.reg} and Proposition~\ref{prop.imp} then show that $\phi$ has the desired regularity,
and statement (A) is proven.

We will now show that statement (B) follows from
statement (A).
If we have a solution as in (A), then we set $g_1:=f^{2/(n-1)}g_0$ with
$f=\<\phi,\phi\>$. Note that $\vol(M,g_1)=\int |\phi|^{2n/(n-1)}=1$.

The transformation formula for the Dirac operator under conformal changes
(Proposition~\ref{prop.conf.change}) implies that there is a
spinor $\phi_1$ on $(M,g_1,\si)$ such that
  $$D_{g_1}\phi_1=\la \phi_1,\qquad |\phi_1|_{g_1}\equiv 1.$$
Then obviously, $\la^+_1(g_1) = \lammin$ and (B) follows.

Statement (C) of Theorem~\ref{theo.main} will follow from 
Proposition~\ref{prop.zero.est}, 
proven in the next section.


\section{The nodal set}\label{sec.degen}

In this section we want to study that the zero set of solutions 
of the Euler-Lagrange equation
\begin{equation}\label{eq.nonlin.p.nn}
  D\phi=c|\phi|^{p-2}\phi.
\end{equation}
Note that very often the zero set of a function $\phi$
solving an equation of the above type is denoted as the 
\emph{nodal set of $\phi$}.

The following theorem is due to C. B\"ar \cite{baer:97} for smooth
$P$.

\begin{theorem}[Nodal sets for Dirac Operators]
Let $(U,g)$ be a Riemannian manifold and let $\phi$ be a solution of
  $$D\phi= P \cdot\phi$$
where $P$ is a smooth function.
Then the nodal set of $\phi$ has Hausdorff dimension at most $n-2$.
\end{theorem}

Unfortunately, \eref{eq.nonlin.p.nn} has not the desired form
as $P=c|\phi|^{p-2}$ is not smooth for non-integer $p$, and
B\"ar's proof does not extend to such a $P$.
Nevertheless, we conjecture that the theorem also holds for 
$P=c|\phi|^{p-2}$, $p>2$. 
\begin{conjecture}
The nodal set of any solution of 
\eref{eq.nonlin.p.nn} has Hausdorff dimension
at most $n-2$.
\end{conjecture}

%
%

If $n=2$ and $p=\pD=4$, then 
we have better regularity. 
In this case solutions of equation~\eref{eq.nonlin.p.nn}
and the corresponding $P=c|\phi|^2$ are smooth.
Hence, using \cite[Main Theorem]{baer:97}, 
one sees that in this case 
the nodal set of a solution is a discrete subset. The following proposition 
controls its cardinality.

\begin{proposition}\label{prop.zero.est}
On a compact spin surface $(M,g,\si)$ of genus $\gamma$ 
let $\phi$ be a solution of equation
  $$D\phi= \lambda |\phi|^2\phi,\quad \|\phi\|_{L^4}=1.$$
Then the number of zeros of $\phi$ is at most
$\gamma-1 + {\lambda^2\over 4\pi}.$
\end{proposition}

In particular, this implies part (C) of Theorem~\ref{theo.main}.
\proof{}
We set $g_1:= |\phi|^4 g$. Then outside the zero set we know by 
Lemma~\ref{lem.scal.bound} that
the Gauss curvature of $g_1$ is at most $\lambda^2$. Furthermore 
$\vol(M,g_1)=1$.
Let $\phi(p)=0$. The integral 
of the geodesic curvature with respect to $g_1$
over small simply closed loop around $p$ is close to $-2 (2j_p+1) \pi$, 
where $j_p$ is the order of the first non-vanishing term in the Taylor 
expansion of $\phi$ in $p$.
We remove small open disks
around the zeros of $\phi$ from $M$, and we obtain a surface with boundary
$M'$. 
With the Gauss-Bonnet theorem we obtain
  $$2\pi\, \chi(M')=\int_{M'} K_{g_1} + \int_{\pa M'} k_{g_1}
           \leq \lambda^2-\sum (2j_p+1) 2\pi.$$
And hence $2\pi\, (2-2\gamma)=2\pi\,\chi(M)\leq\lambda^2- 4\pi \sum j_p$, 
which implies the proposition.
\qed

\section{The spinorial Weierstrass representation}\label{sec.weierrep}

The aim of this section is to recall the spinorial
Weierstrass representation.

Weierstrass published a representation of minimal surfaces
in $\mR^3$ in terms of holomorphic functions 
\cite{weierstrass:66}.
His article deals
only with local questions, everything is described in a fixed
conformal chart of the surface. From a modern (chart free)
point of view, it is clear that
these holomorphic sections should be interpreted as a section of
the spinor bundle, and ``holomorphy'' translates into a ``harmonicity'',
i.e.\ the surface is minimal iff the corresponding spinor $\phi$ satisfies
$D\phi=0$.

During the 20th century several attempts were undertaken to 
globalize the Weierstrass representation and to adapt it to arbitrary surfaces.
Unfortunately, most approaches replaced Weierstrass' original 
approach by a formulation in terms of 
a holomorphic $1$-form and a holomorphic function.
The corresponding formula were quite involved, 
and hence not very suitable for applications. 

An amazing breakthrough was achieved by work of
D.~Sullivan,  R.~Kusner, and N.~Schmitt around 1990, 
and independently by U.~Abresch. 
In early 1989, Dennis Sullivan put together some unpublished notes explaining 
the spinorial character of the Weierstrass representation for minimal 
surfaces. In spring 1989, Robert Kusner realized that the spinor
formalism is not limited to minimal surfaces, but extends to conformal
immersions of arbitrary surfaces, as described below. These techniques were
presented in Sullivan's CUNY seminar in 1992, and 
spread around among the experts rapidly. 
Kusner's results found their continuation
in the PhD thesis of Nick Schmitt \cite{schmitt:diss}.
Schmitt found many interesting applications of the spinorial Weierstrass.
The results of Kusner and Schmitt led to the publications \cite{schmitt:diss}
and \cite{kusner.schmitt:p96}. 

Independently, Abresch developed a spinorial Weierstrass for constant mean
curvature surfaces.
Unfortunately, 
the resulting document of Abresch, some handwritten lecture notes from a 
conference in Luminy, were never published.

We also want to mention an earlier result of Pinkall \cite{pinkall:85a}. 
In the special class of oriented surfaces, Pinkall's result 
establishes a bijection
between regular homotopy classes of immersions of 
oriented surfaces $M$ and $\mZ_2$-valued quadratic forms on $H^1(M,\mZ_2)$.
These quadratic forms are
in bijection with the spin structures obtained in 
the spinorial Weierstrass representation.


\komment{
The idea to use spinors for this representation seems natural 
when one reads Weierstrass' original article 
\cite{weierstrass:66}, but it is very amazing that
the representation takes such a simple formula in terms of spinors. 
Unfortunately, the center of 
mathematicians' interest moved away from Weierstrass' original 
approach to a formulation in terms of 
a holomorphic $1$-form and a holomorphic function, which 
dominated the literature until the late 1980s.

Then, R. Kusner and D. Sullivan found out 
that the spinor formulation gives a simple equation for
non-minimal surfaces, and these ideas led to the preprint 
\cite{kusner.schmitt:p96}.
Although this preprint was not 
published until now, this manuscript is the first document, 
where the spinorial Weierstrass representation was explained.

It also has to be mentioned here that 
independently form Kusner, Sullivan 
and Schmitt, Abresch \cite{abresch:89} had already 
worked with another, less simple version of the spinorial Weierstrass 
representation before, but unfortunately, 
the resulting document, some handwritten lecture notes from a conference in
Luminy, seem to be unavailable to the public.
}

More recent literature concerning this 
representation can be found for example in 
\cite{baer:98,friedrich:98,ammann:diss} and in articles by Pinkall, Taimanov, 
M.U. Schmidt, Morel, Voss and their collaborators. However, this list 
is far from being exhaustive.
 
In our exposition we roughly follow \cite{kusner.schmitt:p96,friedrich:98}. 
The setting for the spinorial Weierstrass representation is as follows.
Let $M$ be a compact Riemann surface of genus $\gamma$. 
The vector bundles 
$\Lambda^{1,0}T^* M$ and $\Lambda^{0,1}T^* M$ are defined 
as the complex linear part and the complex anti-linear part of 
$T^*M\otimes_\mR \mC$. The compositions 
  $$I^{1,0}:T^*M\to T^*M\otimes_\mR \mC\to \Lambda^{1,0}M$$
  $$I^{0,1}:T^*M\to T^*M\otimes_\mR \mC\to \Lambda^{0,1}M$$
of the complexification and the projection on $\Lambda^{1,0} M$
resp.\ $\Lambda^{0,1} M$
define vector space homomorphisms $T^*M\cong \Lambda^{1,0}M$ and  
$T^*M\cong \Lambda^{0,1}M$. These map $I^{1,0}$and $I^{0,1}$ 
preserve the natural 
connections. However, one should pay attention 
to the fact, that these maps do not preserve lengths, but
$2\Ree g(I^{1,0}(\al),I^{1,0}(\be))=2\Ree g(I^{0,1}(\al),I^{0,1}(\be))=g(\al,\be)$. 
The maps $2(I^{1,0})^{-1}$ resp.\ $2(I^{0,1})^{-1}$ is denoted 
as the \emph{real part}, namely $\Ree (\al)=2(I^{1,0})^{-1}(\al)$ for 
$\al\in\Ga(\Lambda^{1,0}M)$ and the same notation is used in the $(0,1)$ case.
\emph{Complex conjugation} maps $\Lambda^{0,1}M$ to $\Lambda^{1,0}M$ 
and vice versa.

As the second Stiefel-Whitney-class of $M$ is the $\mod 2$ reduction 
of the Euler-class of $TM\to M$, one sees that
$M$ is spin. However, the space of spin structures on $M$ is not unique:
it is an affine space 
for the group $H^1(M,\mZ_2)=(\mZ_2)^{2\gamma}$, where $\gamma$ denotes 
the genus of $M$. Hence, there are $4^\gamma$ spin structures on $M$.

If a spin structure is fixed, then the associated vector bundle with respect
to the standard representation of $S^1$ on $\mC$ is a complex 
line bundle $\Si^+M$ satisfying $\Si^+M \otimes \Si^+M\cong \Lambda^{0,1}M$.
If we equip $\Si^+M$ with the natural hermitian metric, we can choose this map such that
the hermitian metric (res.\ the connection) on $\Lambda^{0,1}M$ is 
the tensor product metric (resp.\ tensor product connection).

We define $\Si^-M$ to be $\Si^+M$ with the conjugated complex structure. 
In particular,
there is a natural anti-linear conjugation map $\Si^+M\to \Si^-M$, the 
hermitian product defines a complex bilinear metric contraction 
$\Si^-M\otimes \Si^+M\to \mC$, and $\Si^-M\otimes \Si^-M\cong \La^{1,0}M$. 
In particular,
$\Si^-M\otimes \Lambda^{0,1}M=\Si^-M\otimes \Si^+M \otimes \Si^+M= \Si^+M$, and hence
the Dolbeault operator is a map $\overline\pa:\Gamma(\Si^-M)\to \Gamma(\Si^+M)$, 
and similarly $\overline\pa^*=-\partial: \Gamma(\Si^+M)\to \Gamma(\Si^-M)$.

We define $c^{1,0}$ and $c^{0,1}$ as the compositions
  $$TM\stackrel{b}{\to}T^*M\stackrel{I^{1,0}}\to \Lambda^{1,0}M=\Si^-M\otimes \Si^-M=\Hom_\mC(\Si^+M,\Si^-M),$$
  $$TM\stackrel{b}{\to}T^*M\stackrel{I^{0,1}}\to \Lambda^{0,1}M=\Si^+M\otimes \Si^+M=\Hom_\mC(\Si^-M,\Si^+M).$$
Composing $c^{1,0}$ with complex conjugation yields $c^{1,0}$ and vice versa.
One calculates
  $$c^{1,0}(X) c^{0,1}(Y)+c^{1,0}(Y) c^{0,1}(X)=g(\overline{I(X)},I(Y))+ g(\overline{I(Y)},I(X))=g(X,Y).$$
As a consequence, the map
  $$TM\to\End\left(\Si^+M\oplus \Si^-M\right), \qquad X\mapsto 
\sqrt{2} \begin{pmatrix}
0 &c^{0,1}(X)\\-c^{1,0}(X) & 0 
\end{pmatrix}$$
satisfies the Clifford relations.
\longver{
$$c^{0,1}(X)={1\over 2}\,(X^b- i J(X)^b)$$
$$c^{1,0}(X)={1\over 2}\,(X^b+ i J(X)^b)$$
Assume $e_1=\pa_x$, $e_2=J(e_1)=\pa_y$.
We calculate on $\Si^-M$
  $$e_1\pa_{e_1}= {1\over \sqrt{2}}(dx-i\,dy){\pa\over \pa x}$$
  $$e_2\pa_{e_2}= {1\over \sqrt{2}}(dy+i\,dx){\pa\over \pa y}= {1\over \sqrt{2}}i(dx-i\,dy){\pa\over \pa y}$$
  $$D = \sqrt{2} \,{dx- i\, dy\over 2}\left({\pa\over \pa x}+ i\,{\pa\over \pa y}\right)= 
    \sqrt{2} \,\pa_{\ol z} \,d{\ol z}=\sqrt{2}\,\ol{\pa}$$
}

One sees, that 
the sum $\Si M:=\Sigma^+M\oplus \Sigma^-M$ can be identified with the standard spinor bundle on $M$ in 
such a way that the above map is the Clifford multiplication, and such that
$\Si^+M$ resp.\ $\Si^-M$ are the positive resp.\ negative half-spinors. 

The Dirac operator can be written in this notation as
$$D= \sqrt{2}
\begin{pmatrix}
0 &\overline\partial \\-\partial & 0 
\end{pmatrix}: 
    \Gamma\Big(\Sigma^+\oplus \Sigma^-\Big)\to 
    \Gamma\Big(\Sigma^+\oplus \Sigma^-\Big).$$

Now let us assume that 
$\phi=(\phi_+,\phi_-)$ is a solution of~$D\phi=H|\phi|^2\phi$, where
$H$ is a real-valued function on $M$. 
This means
  $$-\sqrt{2}\pa \phi_+ = H |\phi|^2\phi_-$$
  $$\sqrt{2}\pa \ol{\phi_-} = H |\phi|^2\ol{\phi_+}$$
We define

$$ \alpha := 
{\sqrt{2}}\,
\begin{pmatrix}
  \phi_+\otimes\phi_+ \,+\,\overline{\phi_-}\otimes \overline{\phi_-}\cr
  i \phi_+\otimes \phi_+\,-i\,\overline{\phi_-}\otimes \overline{\phi_-}\cr
  2 i\phi_+\otimes \overline{\phi_-}
\end{pmatrix}\in\Gamma(\Lambda^{0,1}\otimes_\mR \mR^3).
$$
\longver{
$$\pa \al= 2 i H |\phi|^2 
\begin{pmatrix}
-2\Imm (\phi_+\otimes \phi_-) \cr
-2\Ree (\phi_+\otimes \phi_-) \cr
|\phi|^2
\end{pmatrix}\in\Gamma(i\mR^3).$$
}

Let $\witi M$ denote the universal covering of $M$, and $\pi_1(M)$ the group
of Deck transformations.

As $\pa\al=d\al$ is imaginary, we can find a function $F:\witi M\to \mR^3$, 
such that $dF=\Ree \alpha$, and there is a homomorphism
${V:\pi_1(M)\to \mR^3}$, such that 
  $$F(p\cdot\gamma)= F(p) + V(\gamma)\quad \forall p\in \witi M,
    \;\gamma \in \pi_1(M)$$

One calculates that
$F$ is a conformal map with possible branching points,
  $$|dF|=|\Ree \alpha|={1\over \sqrt{2}}\,|\alpha|=|\phi|^2,$$
and that $F(M)$ has mean curvature $H$. 

Hence, the map $F$ satisfies Properties (1) to (3) from the introduction,
i.e. it is a \emph{periodic branched
conformal immersion $F$ based on $(M,g)$} with mean curvature $H$.

In any zero of of the spinor $\phi$, the map $F$ has a branching point.
If $F$ vanishes of order $k$, then $\al$ vanishes of order $2k$.
Hence, all branching points of $F$ are necessarily of even order.

Summarizing the above statement, we obtain for any solution of 
$D\phi=H\,|\phi|^2\phi$  
a periodic branched conformal immersion of $\witi M$ into $\mR^3$ 
which is uniquely 
determined up to translation. If $\phi$ solves $D\phi=H\,|\phi|^2\phi$, 
then  $-\phi$ as well, and the corresponding $F$ is the same. Hence, we obtain
a well-defined map 
$$\left\{\vbox{\hbox{solutions of}\vtop{\hbox{$D\phi=H\;|\phi|^2\phi$}%
    \hbox{on $M$}}}\right\}/{\pm 1} 
\quad\longrightarrow\quad
\left\{\vbox{\hbox{conformal periodic $H$-immersions}\vtop{\hbox{
   $\witi M\to\mR^3$ 
   with branching} 
   \hbox{points of even order}}}\right\}/\vtop{\mbox{trans-}\\\mbox{lations}}$$
and one can show that this map is even a bijection.

The inverse of this map is given by restricting a parallel spinor
on $\mR^3$ to $F(M)$ and
by performing a conformal change \cite{baer:98}.

If $H$ is constant, then 
there is also another version of the spinorial Weierstrass
representation, where the target
space is $S^3$ instead of $\mR^3$. We view $S^3$ as $\SU(2)$ 
with a bi-invariant metric of constant curvature $1$, 
the multiplication in $\SU(2)$ is denoted with 
$\bullet$. The periodicity condition (1) has to be replaced
by
\begin{enumerate}[{\rm (1')}]
\item {\it Left periodicity:}
There is a homomorphism 
$h:\pi(M)\to \SU(2)$, the \emph{periodicity map}, such that 
for any $\gamma\in \pi_1(M)$, and $x\in\witi M$ one has
  $$F(x\cdot \gamma)=h(\gamma)\bullet F(x).$$
Here $\cdot$ denotes the action of $\pi_1$ on $\witi M$ via Deck transformation.\end{enumerate}

One obtains a bijection 
$$\left\{\vbox{\hbox{solutions of}\vtop{\hbox{$D\phi=c\;|\phi|^2\phi$}%
    \hbox{on $M$}}}\right\}/{\pm 1} 
\quad\longrightarrow\quad
\left\{\vbox{\hbox{conformal left periodic $H$-im-}\vtop{\hbox{
mersions $\witi M\to\SU(2)$ 
with} 
   \hbox{branching points of even order}}}\right\}/\vtop{\mbox{Left mul-}\\\mbox{tiplication}}$$
where $c=\sqrt{H^2+1}$.
For details see 
\cite{voss:dipl,morel:02,ammann:habil}.

\section{Applications to constant mean curvature surfaces}\label{sec.app.cmc}

If the dimension of $M$ is $2$, then \eref{eq.dirac.nonlin} reads as 
\begin{equation}\label{eq.nonlin.2}
\Dir\phi=\lammin\, |\phi|^{2} \phi,\qquad \phi\in H_1^{4/3}, \qquad 
\|\phi\|_{L^4}=1
\end{equation}
and according to the regularity theory $\phi$ is even smooth.

The spinorial Weierstrass representation explained in the previous section
tells us, that such a solution can be used to construct 
certain immersions with constant mean curvature.

Combining the previous results we obtain the following application 
that is a stronger version of the 
``Principle for construction of cmc-surfaces''
mentioned in the introduction. 

\begin{proposition}
Assume that the Riemann spin surface $(M,g,\si)$ satisfies 
\begin{equation}\label{ineq.surface}
  \lammin(M,[g],\si)<2\pi.
\end{equation}
Then there is a periodic branched
conformal cmc immersion $F:\witi M\to \mR^3$ based on $(M,g)$. 
The mean curvature is equal to $\lammin(M,[g],\si)$ and the area of a 
fundamental domain is $1$.
The regular homotopy class of $F$ is determined by the spin-structure $\si$.
The indices of all branching points are even, 
and the sum of these indices
is smaller than $2 {\rm genus}(M)$. In particular, 
if $M$ is a torus, there are no branching points.  
\end{proposition}

The proof is a direct consequence of 
Theorem~\ref{theo.main} and the previous section.

There are many examples of stationary points of the functional.
However, it is still open whether they are the maximizers or not.
Note, that by changing the orientation of a surfaces, 
a maximizers $\psi$ of $\cF_\qD$
on $M$ turns into a minimizer $\psi'$ of $\cF_\qD$ on the surface with 
reversed orientation $M'$, with $\cF^{M'}_\qD(\psi')=-\cF^{M}_\qD(\psi)$. 
Let us study some examples.

\begin{examples}
\begin{enumerate}[{\rm (a)}]
\item 
Let $(M,g)$ be a $2$-dimensional torus. Via a conformal change 
we can achieve that $g$ is flat, i.e.\ $M=\mR^2/\Gamma$, equipped 
with the Euclidean metric. 
We assume that the lattice $\Gamma$ is generated by 
$\begin{pmatrix}1\cr 0\end{pmatrix}$ and 
$\begin{pmatrix}x\cr y\end{pmatrix}$, with $y>0$.
The spinor bundle of a flat manifold is flat as well, hence the holonomy
is a map $\Gamma\to \SU(\Si_p M)$. Indeed, the image of this map is contained
in $\{\pm \Id\}$.
We obtain a homomorphism $\chi:\Gamma\to \{\pm \Id\}$.
This homomorphisms
characterizes the spin structure $\si$ in the sense that two spin
structures on $(M,g)$ are isomorphic iff the homomorphisms $\chi$
coincide, and to each such  homomorphism there is a spin structure.
The case $\chi\equiv+\Id$ corresponds to the so-called \emph{trivial} 
spin structure $\si_{\mathrm{tr}}$, the other cases correspond to 
\emph{non-trivial} 
spin structures. 
\footnote{Note that this notation is a bit misleading, as it is only the 
\emph{trivial} spin structure on $M$ that defines a \emph{non-trivial}
element in the bordism class.}

At first, we deal with the case $\si=\si_{\mathrm{tr}}$. 
In this case, after a possible rotation and a possible homothety, 
we can achieve
  $$|x|\leq {1\over 2},\quad 
y^2+x^2\geq 1,\quad y> 0.$$
On $(M,g,\si_{\mathrm{tr}})$ the kernel of $\Dir$ has complex dimension $2$,
and consists of parallel spinors. If one carries out the constructions
from the last sections for a parallel spinor, then one obtains 
an affine conformal map $F$. Such an $F$ is trivially
a periodic conformal immersion $F$ based on $(M,g)$ 
with \emph{vanishing} mean curvature. 
However, if $y>\pi$ then one easily sees that 
$\lammin(M,g,\si_{\mathrm{tr}})<2\pi$. The proposition
yields the existence of a  periodic conformal immersion $F$ 
based on $(M,g)$ with constant mean curvature $\lammin(M,g,\si_{\mathrm{tr}})$,
and the area of a fundamental domain is $1$.
However, the proposition only provides the existence of the solution, but 
we cannot characterize the maximizer.
In our example, a family of solutions can be explicitly written, namely
it is the family of parameterized cylinders 
\begin{eqnarray*}
  F:\mR^2&\to &\mR^3\\
\begin{pmatrix} a\cr b\end{pmatrix}& \mapsto &
P\begin{pmatrix}{\sqrt{y}\over 4\pi}\cos{4\pi b\over y}\cr{\sqrt{y}\over 4\pi}\sin{4\pi b\over y}\cr {a\over \sqrt{y}} \end{pmatrix}+ X_0 
\end{eqnarray*}
for any $P\in O(3)$, $X_0\in \mR^3$. We conjecture that these solutions are
exactly those that correspond to the maximizers and the minimizers of 
$\cF_\qD$, when we normalize such that all spinors have $L^4$-norm $1$.

In the case  $\si\not=\si_{\mathrm{tr}}$
we can achieve
that 
$$\chi\begin{pmatrix}1\cr 0\end{pmatrix}=\Id\qquad \chi\begin{pmatrix}1\cr 0\end{pmatrix}=-\Id,$$
  $$|x|\leq {1\over 2},\quad 
y^2+\left(|x|-{1\over 2}\right)^2\geq {1\over 4},\quad y> 0.$$
The Dirac operator is always invertible.

One easily sees $\lammin(M,g,\si)\leq {\pi\over \sqrt{y}}$.
Hence, the proposition yields solutions for $y>{4\over \pi}$. 
Once again, solutions can be explicitly written, namely
the parameterized cylinder 
\begin{eqnarray*}
  F:\mR^2&\to &\mR^3\\
\begin{pmatrix} a\cr b\end{pmatrix}& \mapsto &
P\begin{pmatrix}{\sqrt{y}\over 2\pi}\cos{2\pi b\over y}\cr{\sqrt{y}\over 2\pi}\sin{2\pi b\over y}\cr {a\over \sqrt{y}} \end{pmatrix}+ X_0 
\end{eqnarray*}
for any $P\in O(3)$, $X_0\in \mR^3$. 
However, in some cases, e.g.\ if $x=0$ and $4/\pi<y<1$, 
these solutions no longer correspond to maximizers and minimizers,
but to saddle points of the functional. We conjecture
that in the case $x=0$, $y<1$ the maximizers and minimizers 
correspond to the unduloid immersions (see Figure~\ref{fig.undu}). 
An unduloid is a surface of revolution of constant mean curvature. 
\begin{figure}
\begin{center}
\bild{
\includegraphics[scale=1]{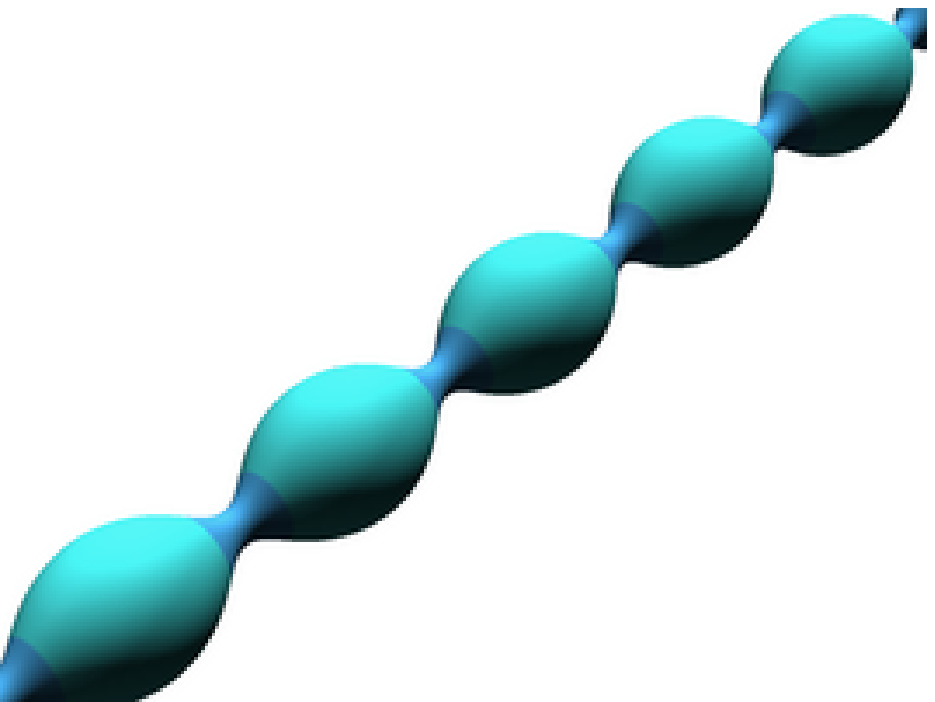}
}
\end{center}
\caption{An unduloid in $\mR^3$, visualized by Nick Schmitt.}\label{fig.undu}
\end{figure}
\item If $M$ has genus $2$, then as in the case of the torus,
the dimension of the kernel is independent of the metric, however it depends
on the spin structure.
If $\si$ is a spin structure such that
$(M,\si)$ is spin-cobordant $0$, then the Dirac operator
is invertible for any metric. Again, as in the torus case, one can find 
for any $\ep>0$ a conformal classes
${g}$ on $M$ with 
$\lammin(M,[g],\si)<\ep$. \cite{ammann.humbert:gluing}
\item If the genus is larger than $2$, then the kernel of the Dirac operator
on a Riemannian spin manifold~$(M,\si)$ 
depends on the metric. For example if $M$ is 
a surface of genus $3$ equipped with the spin structure $\si$ and the conformal
structure ${g_0}$ associated to the periodic conformal immersion with 
vanishing mean curvature indicated in Figure~\ref{fig.kgb.min}.
\footnote{N. Schmitt's illustrations in this paragraph are available 
on his website \url{http://www.gang.umass.edu/gallery/cmc/cmcgallery0302.html},
K. Grosse-Brauckmann's illustrations are available on 
\url{http://www.math.uni-bonn.de/people/kgb/Research/folie\_{}iwp.gif}}
This immersion induces a harmonic spinor on $(M,g,\si)$. However, as 
$(M,\si)$ is spin-cobordant $0$, there is a perturbation $[g_t]$ 
of the conformal structure such that the Dirac operator on $(M,g_t,\si)$
has a trivial kernel for small $t\neq 0$ \cite{maier:97}. In this case 
  $$\lim_{t\to 0\atop t\neq 0}\lammin(M,[g_t],\si)=0,$$ 
hence there exist 
solutions of~\eref{eq.nonlin.2}. Such a solution is visualized in 
Figure~\ref{fig.kgb.cmc}.
\begin{figure}
\begin{center}
\bild{\includegraphics[scale=.5]{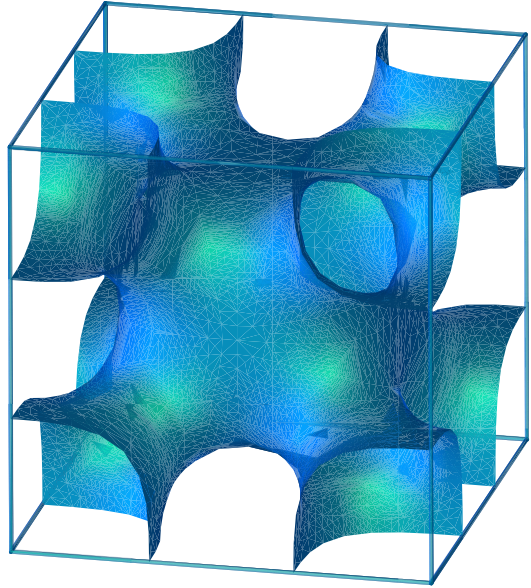}}
\end{center}
\caption{A periodic branched conformal minimal surface%
, visualized by K.~Grosse-Brauckmann}\label{fig.kgb.min}
\end{figure}

\begin{figure}
\begin{center}
\bild{\includegraphics[scale=.5]{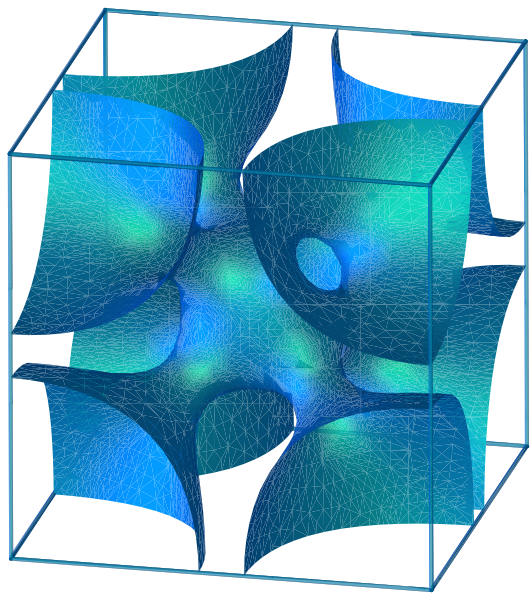}}
\end{center}
\caption{A periodic branched conformal cmc surface%
, visualized by K.~Grosse-Brauckmann}\label{fig.kgb.cmc}
\end{figure}
\item Constant  mean curvature immersions of $T^2$ into $\mR^3$, 
in particular Wente tori and twisty tori (Figure~\ref{fig.twisty}) 
also correspond to 
stationary points of $\cF_\qD$. However, as they are not embedded  
\cite{li.yau:82} tells us that $\int H^2 \geq 8\pi$.
On the other hand maximizers and minimizers of $\cF_\qD$ satisfy 
$\int H^2= \lammin^2\leq 4\pi$, hence these tori do \emph{not} 
correspond to maximizers or minimizers.   
\begin{figure}
\begin{center}
\bild{
\includegraphics[scale=1]{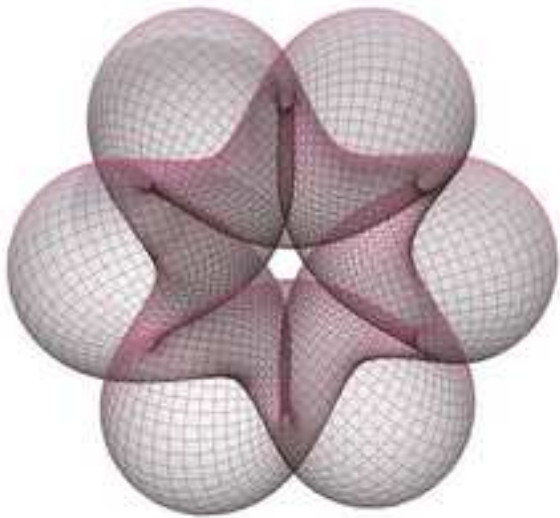}
}\end{center}
\caption{A twisty torus, a cmc immersed torus in $\mR^3$,
visualized by Nick Schmitt}\label{fig.twisty}
\end{figure}
\end{enumerate}
\end{examples}

Similar propositions also hold for immersions into $S^3$ and 
into hyperbolic space $H^3$, see \cite{ammann:habil}. We will only 
specify the case of $S^3=\SU(2)$.

\begin{proposition}
Assume that the Riemann spin surface $(M,g,\si)$ satisfies 
\begin{equation}
  \lammin(M,[g],\si)<2\pi.
\end{equation}
Let $a\in(0,\lammin(M,[g],\si)$ be given.
Then there is a left periodic branched
conformal cmc immersion $F:\witi M\to \SU(2)$ based on $(M,g)$. 
The mean curvature is equal to $H=\sqrt{(\lammin(M,[g],\si)/a)^2-1}$, 
and the area of a 
fundamental domain is $a^2$.
The regular homotopy class of $F$ is determined by the spin-structure $\si$.
The indices of all branching points are even, 
and the sum of these indices
is smaller than $2 {\rm genus}(M)$. In particular, 
if $M$ is a torus, there are no branching points.  
\end{proposition}

\begin{figure}
\begin{center}
\bild{
\includegraphics[scale=.3]{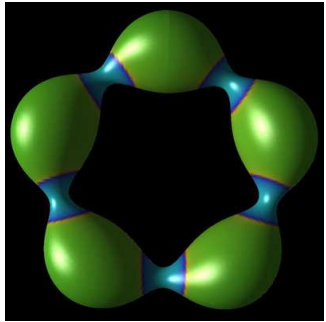}
}\end{center}
\caption{An unduloid in $S^3$, visualized by Nick Schmitt.}\label{fig.undu.s}
\end{figure}

An example where the image of the periodicity map has a finite image in 
$\SU(2)$ is given in Figure~\ref{fig.undu.s}.

\begin{remark}
If $(M,g)$ is an analytic Riemannian manifold of dimension $3$, then there
is an analogue of the spinorial Weierstrass representation. This implies 
that solutions $\phi$ of \eref{eq.nonlin.p} 
can be geometrically interpreted as 
a conformal cmc embedding of $M\setminus \phi^{-1}(0)$ 
into a (non-complete) $4$-dimensional Riemannian manifold $(N,h)$.
This manifold $(N,h)$ carries a parallel spinor whose restriction
to $(M,g)$ is again the solution of~\eref{eq.nonlin.p}. The manifold 
$(N,h)$ depends on $\phi$ and is unique up to restriction to subsets 
and coverings. A more detailed exposition of this $3$-dimensional version
is work in progress.
\end{remark}

\begin{appendix}

\section{Elliptic regularity}\label{sec.regularity}

In this section we want to collect some facts about elliptic regularity
for Dirac operators. The proofs of these statements are
analogous to proofs of the corresponding statements for the Laplace operators
as done e.g.\ in \cite{gilbarg.trudinger:77} and \cite{adams:75}.
Details on how to prove most statements of this section
are given also in \cite{ammann:habil}. Either a proof is provided there or 
it is sketched how the statements are reduced to standard theorems.

Let $(N,g)$ be a Riemannian manifold, possibly with boundary,
with a spin structure~$\si$. The interior of $N$ is denoted as $N_0$.
The spinor bundle $\Si(N,g,\si)$ is a complex vector bundle
carrying a natural connection and a natural hermitian metric.

We now define the Sobolev norms of spinors fields.

\begin{definition}[Sobolev spaces]
For any spinor $\psi$ smooth on $N$, 
and $q\in (1,\infty)$, $k\in \mN\cup\{0\}$,
we define the $H^q_k$-norm of $\psi$ as
  $$\|\psi\|_{H^q_k(N,g,\si)}
   = \sum_{l=0}^k\,\Biggl(\int_N\Bigl|\underbrace{\na \ldots \na}_{l-\mbox{times}}\psi\Bigr|^q\,\dvol_g\Biggr)^{1/q}.$$
When the domain of integration is clear, we simply write 
$\|\psi\|_{H^q_k}$ instead of $\|\psi\|_{H^q_k(N,g,\si)}$.
The closure of the space of 
smooth spinors with respect to this norms is denoted as $H^q_k(\Si(N,g,\si))$.
If the underlying structure of Riemannian spin manifold
is clear from the context, we often write shortly 
$H^q_k$,  or in other situations where we want to emphasize the metric
we write  $H^q_k(g)$. We will also write 
$L^q$ for $H^q_0$. 
\end{definition}

\begin{theorema}[Interior $L^p$ estimates]\label{theo.int.lp}
Let $(N,g,\si)$ be a compact Riemannian spin 
manifold, possibly \emph{with boundary}. Let $K$ be a compact
subset of $N_0$.
Let $\psi$ be in $H_k^q(\Si(N,g,\si))$ and let $\phi\in H_1^1(\Si(N,g,\si))$ 
be a weak solution (i.e.\ in the sense of distributions) of
  $$D\phi = \psi$$
on $N_0$.
Then $\phi|_K\in H_{k+1}^q(\Si(K,g,\si))$ and
  $$\|\phi\|_{H_{k+1}^q(\Si(K,g,\si))}\leq C\cdot
    \left(\|\psi\|_{H_k^q(\Si(N,g,\si))} + \|\phi\|_{L^q(\Si(N,g,\si))}\right),$$
where $C=C(N,g,\si,K,k,q)$, i.e.\ $C$ depends only on $N,g,\si,K,k$ and $q$.

Furthermore, 
if a sequence of metrics $(g_i)_{i\in \mN}$ 
converges in the $C^\infty$-topology
to a Riemannian metric $g$, then the constants $C$ can be chosen 
such that 
  $$\sup_{i\in \mN}  C(N,g_i,\si,K,k)<\infty.$$   
\end{theorema}

In order to prove the theorem, it is sufficient to prove 
it in the case that $N$ and $K$ are small concentric geodesic balls, the radius
being controlled in terms of curvature bounds, and bounds on derivatives
of the curvature. As explained in \cite{ammann:habil} 
this can be done in a way analogous to the corresponding statement 
for the Laplacian in \cite{gilbarg.trudinger:77}.
The general statement then can be reduced to the special 
case by covering  (the general) $K$ by finitely many small open balls and an 
associated partition of unity.

In the special case that the boundary of $N$ is the empty set,
and that $K=N=N_0$, 
then there is a stronger version, that will be used as well.

\begin{theorema}[Global $L^p$ estimates]\label{theo.global.lp}
Let $(N,g,\si)$ be a compact Riemannian spin manifold 
\emph{without boundary} and $\psi\in H_k^q(\Si(N,g,\si))$.
Then any weak solution  $\phi\in H_1^1(\Si(N,g,\si))$ of
  $$D\phi = \psi$$
satisfies $\phi\in H_{k+1}^q(\Si(N,g,\si))$, and there is a constant $C=C(N,g,\si)$
such that
  $$\|\phi\|_{H_{k+1}^q}\leq C 
    \left(\|\psi\|_{H_{k}^q}+\|\pi_{\mathop{\rm ker}(D)}(\phi)\|_{L^q}\right)$$
where $\pi_{\mathop{\rm ker}(D)}$ is the $L^2$-orthogonal projection 
to the kernel of $D$.
\end{theorema}

The proof of this theorem is not difficult if one uses some facts 
about pseudodifferential operators, explained for example in 
\cite{taylor:81}. We abbreviate  $\pi:=\pi_{\mathop{\rm ker}(D)}$.
The spectrum of the elliptic operator  
$D+\pi$ is bounded away from $0$,
hence it is invertible and its inverse is a pseudodifferential operator
of degree $-1$. Pseudodifferential operators of degree $-1$
are continuous from $L^q$ to $H_1^q$. Hence, 
  $$\|\phi\|_{H_1^q}=\|(D+\pi)^{-1}(D+\pi)\phi\|_{H_1^q}\leq C\|(D+\pi)\phi\|_{L^q}\leq C\left(\|D\phi\|_{L^q}+\|\pi(\phi)\|_{L^q}\right).$$
This is the statement of the theorem in the case $k=0$, and the statement
for $k>0$ follows for example by using Theorem~\ref{theo.int.lp}.

However, a uniformity statement for converging $g_i\to g$ as in 
Theorem~\ref{theo.int.lp} does not hold. If 
$\dim \ker D_g >\limsup_i\dim \ker D_{g_i}$, then one easily shows that
there are eigenspinors $\phi_i$ for $D_{g_i}$ to eigenvalues 
$\la_i\to 0$, $\la_i\neq 0$ with  
  $${\|D\phi_i\|_{L^q(g_i)}+\|\pi_{g_i}(\phi_i)\|_{L^q(g_i)}\over \|\phi_i\|_{H_1^q(g_i)}}=\la_i,$$
which would contradict a uniform version of Theorem~\ref{theo.global.lp}.

We also needs some facts about H\"older norms.

\begin{definition}[H\"older spaces]
Let $\al\in (0,1]$.
On $C^\infty(\Si(N,g,\si))$ we define the
\emph{H\"older norms}
\begin{eqnarray*}
  \|\phi\|_{C^{0,\al}} & := & \mbox{h\"ol}_\al(\phi)\\
  \|\phi\|_{C^{1,\al}} & := & \|\phi\|_{C^0} + \mbox{h\"ol}_\al(\nabla \phi)\\
  \mbox{h\"ol}_\al(Q)& := & \sup \Bigl\{{|Q(x)-P_\ga Q(y)|\over
   d(x,y)^\al}\;|\; x,y\in M, x\neq y, P_\ga \mbox{ is the }\\
 && \mbox{parallel transport along a shortest geodesic $\ga$ from $x$ to $y$}.\Bigr\}
\end{eqnarray*}
If the shortest geodesic is not unique, then we use the convention that
the supremum runs over all possible choices of shortest geodesics.

The completions of $C^\infty(\Si(N,g,\si))$ with respect to the 
H\"older norms $C^{0,\al}$ and $C^{1,\al}$ define the 
\emph{H\"older spaces}
$C^{0,\al}=C^{0,\al}(\Si(N,g,\si))$ and
$C^{1,\al}=C^{1,\al}(\Si(N,g,\si))$.
\end{definition}

H\"older norms are important as they admit Schauder estimates. 

\begin{theorema}[Schauder estimates]\label{theo.schauder}
Let $(N,g,\si)$ be a compact Riemannian spin manifold, possibly
\emph{with boundary}, and let $K$ be a compact subset of $N_0$. 
Suppose $\psi\in C^{k,\al}(\Si(N,g,\si))$, $k\in \mN\cup\{0\}$.
Then for any weak solution $\phi\in L^1 (\Si(N,g,\si))$ of
  $$D\phi = \psi$$
we have $\phi|_K\in C^{k+1,\al}(\Si(K,g,\si))$ and
  $$\|\phi\|_{C^{k+1,\al}(\Si(K,g,\si))}\leq 
    C \cdot (\|\psi\|_{C^{k,\al}(\Si(N,g,\si))}
    +\|\phi\|_{C^0(\Si(N,g,\si)})$$
where $C=C(N,g,\si,K,\al)$.

Furthermore, 
if a sequence of metrics $(g_i)_{i\in \mN}$ 
converges in the $C^\infty$-topology
to a Riemannian metric $g$, then the constants $C$ can be chosen 
such that 
  $$\sup_{i\in \mN}  C(N,g_i,\si,K,k)<\infty.$$   
\end{theorema}

As before, one important special case is that $N$ has empty boundary
and $K=N=N_0$.

The proof can be done in a way analogous
to the proof of the corresponding statements for the Laplacian in 
\cite{gilbarg.trudinger:77}. Again it is sufficient to prove it for 
small concentric geodesic balls, and to glue them together.
We omit the details. For a proof for concentric geodesic 
balls is provided for example in \cite{ammann:habil}.

In local charts one easily reduces the following theorem to the standard
Arcela-Ascoli theorem.

\begin{theorema}[Arcela-Ascoli] \label{theo.arc.asc}
Let $(N,g,\si)$ be a Riemannian spin
manifold, possibly \emph{with boundary}. 
For $m\in \mN\cup \{0\}$, the inclusion
$C^{m,\al}(\Si(N,g,\si))\to C^m(\Si(N,g,\si))$ is compact, i.e.\ 
a bounded sequence in $C^{m,\al}(\Si(N,g,\si))$ has a subsequence 
convergent in $C^m(\Si(N,g,\si))$.
\end{theorema}

Sobolev and H\"older spaces are related by a several embedding theorems,
some of them are compact, others are only bounded. We summarize the embeddings
that are needed in the article.

\begin{theorema}[Embedding theorems]\label{theo.sobo}\
Let $k,s\in \mN\cup\{0\}$, $k\geq s$ and $q,r\in (1,\infty)$.
Let $(M,g,\si)$ be a compact Riemannian spin manifold \emph{without boundary}.
All spaces of functions are defined on sections of $\Si(M,g,\si)$.
\begin{enumerate}[{\rm (a)}]
\item {\rm (Sobolev embedding theorem I).} If
 $$\phantom{\mathrm{(A.1)}}\qquad\qquad\qquad\qquad {1\over r} - {s\over n}\geq {1\over q} - {k\over n},\qquad\qquad\qquad\qquad \mathrm{(A.1)}$$
then $H_k^q$ is continuously embedded into $H_s^r$.
\item {\rm (Rellich-Kondrakov theorem).}
If strict inequality holds in {\rm (A.1)} 
and if $k>s$, then the inclusion
$H_k^q\embed H_s^r$ is a
compact map.
\item  {\rm (Sobolev embedding theorem II).}
Suppose $0<\al <1$, $m\in\{0,1\}$ and
  $${1\over q}\leq {k-m-\al \over n}.$$
Then $H_k^q$ is continuously embedded into $C^{m,\al}$.
\end{enumerate}
\end{theorema}

\section{Some facts about H\"older spaces}\label{sec.scha.power}

In this section we want to include some proofs of probably 
well-known statements about H\"older spaces.

\begin{lemmaa}
Let $V$ be a Euclidean vector space, $\al\in (0,1)$, $m\in \mN$. 
Then the map $V\to \bigotimes_m V$, 
  $$T:x\mapsto |x|^{\al-m}\underbrace{x\otimes\cdots\otimes x}_{\mbox{$m$ times}}$$
is $C^\al$.
\end{lemmaa}

\proof{}
Let $x,y\in V$. At first, we suppose that
$\|x-y\|<\de$ and $\|x\|\geq \de$.
Then we calculate
$$
\|x\|\,\Big\|{x\over \|x\|}-{y\over \|y\|}\Big\|=\Big\|x-{\|x\|\over \|y\|}y\Big\|\leq \|x-y\|+\Big|1- {\|x\|\over \|y\|}\Big|\,\|y\|\leq 2\de
$$
Hence 
$$
   \Big\|{x\over \|x\|}\otimes \cdots \otimes {x\over \|x\|} 
    - {y\over \|y\|}\otimes \cdots \otimes {y\over \|y\|}\Big\|\leq {2\de m\over \|x\|}
$$
This implies 
  $$|T(x)-T(y)|\leq 
    \Big|\|x\|^\al -\|y\|^\al\Big|+ \|x\|^\al \,{2\de m\over \|x\|}\leq \|x-y\|^\al+ \|x\|^{\al-1}2\de m\leq (2m+1)\de^\al.$$

Now suppose that
$\|x-y\|<\de$ and $\|x\|< \de$. Then 
  $$|T(x)-T(y)|\leq \|x\|^\al + \|y\|^\al\leq 3\de^\al.$$
Hence, we obtain $|T(x)-T(y)|\leq (2m+1)\|x-y\|$ for all $x,y\in V$, and hence $T$ is $C^\al$.
\qed

\begin{propositiona}
If $V$ is a vector bundle over a compact Riemannian manifold, 
and if $\phi$ is a
$C^{1,\al}$-section of $V$, $\al\in (0,1)$ then 
$\psi=|\phi|^\beta\phi$, $\be>0$ is a $C^{1,\gamma}$-section
for $\gamma:=\min\{\al,\be\}$.
\end{propositiona}

\proof{}
The section $\phi$ is obviously Lipschitz, hence $|\phi|^\be$ is $C^\be$.
We have to show that 
\begin{eqnarray*} 
  \na\psi&=&|\phi|^\be\na\phi + \beta\<\na \phi,\phi\>|\phi|^{\be-2}\phi
\end{eqnarray*}
is $C^\ga$. The first summand is a product of $C^\be$ and $C^\al$, 
hence $C^\gamma$. According to the previous lemma, 
$|\phi|^{\be-2}\phi\otimes \phi$ is $C^\be$, hence the second summand
is a product of $C^\al$ and $C^\be$, hence also $C^\ga$.
\qed

\end{appendix}


\providecommand{\bysame}{\leavevmode\hbox to3em{\hrulefill}\thinspace}
\providecommand{\MR}{\relax\ifhmode\unskip\space\fi MR }
\providecommand{\MRhref}[2]{%
  \href{http://www.ams.org/mathscinet-getitem?mr=#1}{#2}
}
\providecommand{\href}[2]{#2}

\pagebreak[2]

\vspace{1cm}
Author's address:
\nopagebreak
\vspace{5mm}\\
\parskip0ex
\vtop{
\hsize=7cm\noindent
\obeylines
Bernd Ammann
Institut \'Elie Cartan BP 239              
Universit\'e Henri Poincar\'e, Nancy 1               
54506 Vandoeuvre-l\`es-Nancy Cedex               
France                           
}

\vspace{0.5cm}

E-Mail:
{\tt bernd.ammann@gmx.de}

WWW:
{\tt http://www.berndammann.de/uni}

\end{document}